\documentclass[10pt,reqno]{amsart}
 \usepackage{graphicx,enumitem,dsfont}
\usepackage{a4wide}
\usepackage{color}
\usepackage{upref}
\usepackage{stmaryrd}
\usepackage{caption}
\usepackage{subcaption}

\newcommand {\vv}{\llbracket v \rrbracket}

\newcommand{\W}{{\mathcal{W}}}
\newcommand{\n}{{\mathbf{n}}}

 \newcommand{\R}{\mathbb{R}}

 \newcommand{\Comment}[1]{}

 \newcommand*{\jump}[1]{\llbracket #1 \rrbracket}
 \newcommand{\eps} {\varepsilon}
 \renewcommand{\i}{\ifmmode\mathit{\mathchar"7010 }\else\char"10 \fi}
 \renewcommand{\j}{\ifmmode\mathit{\mathchar"7011 }\else\char"11 \fi}

 \newcommand{\test}{\varphi}
 \newcommand{\Dx}{\Delta x}
 \newcommand{\Dy}{\Delta y}

 \newcommand{\abs}[1]{\left|#1\right|}

 \newcommand{\Pt}{\Pi_T}

 \newtheorem{definition}{Definition}[section]
 \newtheorem{theorem}{Theorem}[section]
 \newtheorem*{maincorollary*}{Main Corollary}

 \newtheorem{remark}{Remark}[section]
 
 \newtheorem*{maintheorem*}{Main Theorem}
 \allowdisplaybreaks
 \numberwithin{equation}{section}
 \numberwithin{figure}{section}
 \numberwithin{table}{section}

{%

\begin{enumerate}}%
{\end{enumerate}}

\newenvironment{Definitions}% Definition of definitions
{%

\begin{enumerate}}%
{\end{enumerate}}

\title[RKDG schemes for a variational wave equation]
{High-order energy stable numerical schemes
          for a \\nonlinear variational wave equation 
          modeling \\nematic liquid crystals in two dimensions}

\author[P. Aursand]{Peder Aursand}
\address[Peder Aursand]{\newline
    Department of Mathematical Sciences,
    \newline Norwegian University of Science and Technology,
   \newline  NO--7491 Trondheim, Norway.}
\email[]{peder.aursand@math.ntnu.no}

\author[U. Koley]{Ujjwal Koley} \address[Ujjwal Koley] {\newline  
   Tata Institute of Fundamental Research,  \newline
   Centre For Applicable Mathematics,   
\newline Post Bag No. 6503, GKVK Post Office, \newline Sharada Nagar, Chikkabommasandra, \newline
Bangalore 560065, India.} 
\email[]{ujjwal@math.tifrbng.res.in}

\date{\today}

\begin{document}

\begin{abstract}
  We consider a nonlinear variational wave equation that models the dynamics
  of the director field in nematic liquid crystals with high molecular
  rotational inertia. Being derived from an energy principle, energy stability
  is an intrinsic property of solutions to this model. For the two-dimensional
  case, we design numerical schemes based on the discontinuous Galerkin
  framework that either conserve or dissipate a discrete version of the
  energy. 
  
  Extensive numerical experiments are performed verifying the scheme's energy
  stability, order of convergence and computational efficiency. 
  The numerical solutions are compared to those of a simpler first-order Hamiltonian scheme.
  We provide numerical evidence that solutions of the 2D variational wave
  equation loose regularity in finite time. After that occurs,
  dissipative and conservative schemes appear to converge to different solutions. 
\end{abstract}

\maketitle

%\tableofcontents

\section{Introduction}
\label{sec:intro}
\subsection{The Equation} 
Liquid crystals (LCs) are mesophases, i.e., intermediate states of matter between
the liquid and the crystal phase. They possess some of the
properties of liquids (e.g. formation, fluidity) as well as some crystalline
properties (e.g. electrical, magnetic, etc.) normally associated with solids. 
The nematic phase is the simplest of the liquid crystal mesophases, and is close to 
the liquid phase. It is characterized by long-range orientational order, i.e., the 
long axes of the molecules tend to align along a preferred direction, which can be 
considered invariant under rotation by an angle of $\pi$. 
The state of a nematic liquid crystals is usually given by two linearly independent 
vector fields; one describing the fluid flow and the other describing the dynamics 
of the preferred axis, which is defined by a vector $\n$ giving its local orientation. 
Under the assumption of constant degree of orientation, the magnitude of the 
\emph{director field} $\n$ is usually taken to be unity. In the present work
we focus exclusively on the dynamics of the director field 
(independently of any coupling with the fluid flow), a map
\begin{align*}
  \n : \R^3 \times [0, \infty) \rightarrow \mathbb{S}^2
\end{align*}
from the Euclidean space to the unit ball.

We consider the elastic dynamics of the liquid crystal director field in the
inertia-dominated case (zero viscosity). Associated with the director field $\n$, 
the classical Oseen-Frank elastic energy density $\W$ is given by
\begin{equation}
  \begin{aligned}
    \W(\n, \nabla \n) = \alpha \abs{\n \times
      (\nabla \times \n)}^2 + \beta \left(\nabla \cdot
      \n \right)^2 + \gamma \left( \n \cdot (\nabla
      \times \n) \right)^2.
  \end{aligned}
  \label{eq:osf}
\end{equation}
The constants $\alpha, \beta$ and $\gamma$ are elastic material constants of the 
liquid crystal, and are associated with the three basic types of deformations of the
medium; bend, splay and twist; respectively.  Each of these constants must be
positive in order to guarantee the existence of the minimum configuration of
the energy $\W$ in the undistorted nematic configuration.

The one constant approximation ($\alpha = \beta = \gamma$) often
provides a valuable tool to reach a qualitative insight into distortions of
nematic configurations. Observe that, in this case the potential energy
density \eqref{eq:osf} reduces to the Dirichlet energy
\begin{align*}
  \W(\n, \nabla \n) = \alpha \abs{\nabla
    \n}^2.
\end{align*}
This corresponds to the potential energy density used in harmonic
maps into the sphere $\mathbb{S}^2$. 
The stability of the general
Oseen--Frank potential energy equation, derived from the potential
\eqref{eq:osf} using a variational principle, is studied by 
Ericksen and Kinderlehrer~\cite{ericksen}.
For the parabolic flow associated to \eqref{eq:osf}, see
\cite{Beres,coron} and references therein.

In the regime in which inertial effects dominate viscosity, the dynamics of the
director $\n$ is governed by the least action principle
\begin{equation}
  \begin{aligned}
  \mathbb{J}(\n) =  \iint \left(\n_t^2 -
      \W(\n, \nabla \n) \right) \,dx\,dt, \qquad
    \n \cdot \n =1.
  \end{aligned}
  \label{eq:lap}
\end{equation}
Standard calculations reveal that the \emph{Euler-Lagrange} equation
associated to $\mathbb{J}$ is given by
\begin{equation}
\label{eq:euler}
\begin{aligned} 
\n_{tt} = \mathrm{div} \left(  \W_{\nabla \n} (\n, \nabla \n)\right) -
\W_{\n} (\n, \nabla \n),
\end{aligned}
\end{equation}
and is termed the variational wave equation.
Introducing the \emph{energy} and \emph{energy density}
\begin{align*}
\mathcal{E}(t) = \int \left(  \n_t^2 + \W (\n, \nabla \n)\right) \,dx, \qquad \mathbf{E}(t,x) =  \n_t^2 + \W (\n, \nabla \n),
\end{align*}
it is easy to check the identities 
\begin{align*}
\mathcal{E}' =0, \qquad \mathbf{E}_t = \mathrm{div} \left( \W_{\nabla \n} (\n,
\nabla \n) \n_t \right),
\end{align*}
in light of \eqref{eq:euler}.
Given the formidable difficulties in the mathematical analysis of
\eqref{eq:euler}, it is customary to investigate the particular case of a
planar director field configuration.

The physical implications of considering the inertia-dominated regime warrants
a comment. Indeed, in many experimental situations the inertial forces acting
on the director are orders of magnitude smaller than the dissipative. For this
reason, the inertial term is often neglected in modelling
\cite{Stewart2004,vanDoorn1975,Gang1987}. It was however noted early by
Leslie~\cite{Leslie1979} that inertial forces might be significant in cases
where the director field is subjected to large accelerations. In general,
inertia will be more significant in the small time-scale dynamics of the
director.  For this reason, their inclusion can be warranted in, e.g., liquid
crystal acoustics \cite{Kapustina2004}, mechanical vibrations
\cite{Vladimirov2007} and in cases with and external oscillating magnetic
field \cite{Yun1973}.

\subsubsection{One-dimensional planar waves}
Planar deformations are central in the mathematical study of models for nematic 
liquid crystals. A simple such model can be derived by assuming that the deformation
depends on a single space variable $x$ and that the director field
$\n$ in confined to the $x$-$y$ plane. In this case we can write the director
as
\begin{align*}
  \n= (\cos u(x,t),\, \sin u(x,t), \, 0).
\end{align*}
Geometrically, the molecules are lined up vertically on the $x$-$y$ plane, and at each
column (located at $x$) $u(x,t)$ measures the angle of the director field to the
$x$-direction. 
With the above simplifications, the variational principle \eqref{eq:lap} reduces to
\begin{equation}
  \label{eq:main}
  \begin{cases}
    u_{tt} - c(u)\left( c(u) u_x\right)_x  =0, \quad (x,t) \in \Pt, &  \\
    u(x,0) = u_0(x),\quad x \in \R, &  \\
    u_t(x,0) = u_1(x), \quad x \in \R, &
  \end{cases}
\end{equation}
where $ \Pt = \R \times [0,T]$ with fixed $T>0$ , and the wave speed
$c(u)$ given by
\begin{equation}
  \label{eq:wavespeed}
  c^2(u) = \alpha \cos^2 u + \beta \sin^2u.
\end{equation}
Initially considered by Hunter and Saxton \cite{saxton,hs1991}, 
\eqref{eq:main} is the simplest form of the nonlinear variational wave equation
\eqref{eq:euler} studied in the literature.

\subsubsection{Two-dimensional planar waves}
Planar deformations can also be studied in two dimensions. Specifically, 
if we assume that the deformation depends on two space variables $x,y$,
the director can be written in the form
\begin{align*}
  \n = \left(\cos u(x,y,t), \, \sin u(x,y,t), \, 0 \right)
\end{align*}
with $u$ being the angle to the $x$-$z$ plane. The corresponding variational wave
equation is given by
\begin{equation}
  \label{eq:main_2d}
  \begin{cases}
    u_{tt} - c(u)\left( c(u) u_x\right)_x -b(u)\left( b(u) u_y\right)_y - a'(u) u_x u_y - 2 a(u) u_{xy} =0, \quad (x,y,t) \in \mathbb{Q}_T, &  \\
    u(x,y,0) = u_0(x,y),\quad (x,y) \in \R^2, &  \\
    u_t(x,y,0) = u_1(x,y), \quad (x,y) \in \R^2, &
  \end{cases}
\end{equation}
where $ \mathbb{Q}_T = \R^2 \times [0,T]$ with $T>0$ fixed, $u:
\mathbb{Q}_T \rightarrow \R$ is the unknown function and $a,b,c$ are
given by
\begin{align*}
  c^2(u) & = \alpha \cos^2 u + \beta \sin^2 u, \\
  b^2(u) &= \alpha \sin^2 u + \beta \cos^2 u, \\
  a(u) &= \frac{\alpha - \beta}{2} \sin(2u).
\end{align*} 
In this picture, $c(u)$ is the wave speed in the $x$-direction and $b(u)$ is
the wave speed in the $y$-direction.
%%In particular, we have that $b^2 + c^2 = \alpha + \beta$ and $a'' = -4a$.

For smooth solutions of \eqref{eq:main_2d} it is straightforward to verify that
the energy
\begin{equation}
  \begin{aligned}
    \label{eq:energy_2d}
    \mathcal{E}(t) &= \iint_{\R^2} \left( u_t^2 + c^2(u) u_x^2 + b^2(u) u_y^2 + 2 a(u) u_x u_y \right) \,dx \,dy \\
    &= \iint_{\R^2} u_t^2 + \left( \alpha( \cos(u) u_x + \sin(u)
      u_y) \right)^2 + \left(\beta( \sin(u) u_x - \cos(u) u_y) \right)^2\,dx\,dy
  \end{aligned}
\end{equation}
is conserved, i.e., we have
\begin{equation}
  \label{eq:econ}
  \frac{d \mathcal{E}(t)}{dt} \equiv 0.
\end{equation}
Moreover, for all $t \in [0,T]$ we have 
\begin{align*}
\iint_{\R^2} \left(u_t^2 + \min\lbrace\alpha, \beta \rbrace (u_x^2 + u_y^2) \right) \,dx\,dy \le \mathcal{E}(t) \le 
\iint_{\R^2} \left(u_t^2 + \max \lbrace \alpha, \beta\rbrace (u_x^2 + u_y^2) \right) \,dx\,dy.
\end{align*}
In particular, it follows that $\mathcal{E}(t) \ge 0$ for all $t \in [0,T]$.
To see this, first we consider $\alpha \ge \beta$ (for $\alpha \le \beta$, we argue in the same way). Then
\begin{align*}
 c^2(u) u_x^2 &+ b^2(u) u_y^2 + 2 a(u) u_x u_y \\
&= \left(\alpha \cos^2 (u) + \beta \sin^2 (u)  \right) u_x^2 + \left( \alpha \sin^2 u + \beta \cos^2 u \right) u_y^2 + 2 (\alpha - \beta) \sin (u) \cos(u) u_x u_y \\
& \le \left(\alpha \cos^2( u) + \beta \sin^2 (u)  \right) u_x^2 + \left( \alpha \sin^2 u + \beta \cos^2 u \right) u_y^2 + 2 (\alpha - \beta) \abs{\sin (u) \cos(u) u_x u_y} \\
& \qquad = \alpha \left( \cos^2 (u) u_x^2 + \sin^2 (u)  u_y^2 + 2 \abs{\sin (u) \cos(u) u_x u_y} \right) \\
& \qquad \qquad +\beta \left(  \sin^2 (u) u_x^2 + \cos^2 (u)  u_y^2 - 2 \abs{\sin (u) \cos(u) u_x u_y}  \right) \\
& =\alpha \left( \abs{\cos (u) u_x} + \abs{\sin (u)  u_y}  \right)^2 + \beta \left( \abs{ \sin (u) u_x} -\abs{ \cos (u)  u_y}  \right)^2 \\
& \le \alpha \Big[ \left( \abs{\cos (u) u_x} + \abs{\sin (u)  u_y}  \right)^2 +  \left( \abs{ \sin (u) u_x} -\abs{ \cos (u)  u_y}  \right)^2 \Big] = \alpha (u_x^2 + u_y^2),
\end{align*}
and 
\begin{align*}
 c^2(u) u_x^2 &+ b^2(u) u_y^2 + 2 a(u) u_x u_y \\
&= \left(\alpha \cos^2 (u) + \beta \sin^2 (u)  \right) u_x^2 + \left( \alpha \sin^2 u + \beta \cos^2 u \right) u_y^2 + 2 (\alpha - \beta) \sin (u) \cos(u) u_x u_y \\
& \ge \left(\alpha \cos^2( u) + \beta \sin^2 (u)  \right) u_x^2 + \left( \alpha \sin^2 u + \beta \cos^2 u \right) u_y^2 - 2 (\alpha - \beta) \abs{\sin (u) \cos(u) u_x u_y} \\
& \qquad = \alpha \left( \cos^2 (u) u_x^2 + \sin^2 (u)  u_y^2 - 2 \abs{\sin (u) \cos(u) u_x u_y} \right) \\
& \qquad \qquad +\beta \left(  \sin^2 (u) u_x^2 + \cos^2 (u)  u_y^2 + 2 \abs{\sin (u) \cos(u) u_x u_y}  \right) \\
& =\alpha \left( \abs{\cos (u) u_x} - \abs{\sin (u)  u_y}  \right)^2 + \beta \left( \abs{ \sin (u) u_x} +\abs{ \cos (u)  u_y}  \right)^2 \\
& \ge \beta \Big[ \left( \abs{\cos (u) u_x} - \abs{\sin (u)  u_y}  \right)^2 +  \left( \abs{ \sin (u) u_x} +\abs{ \cos (u)  u_y}  \right)^2 \Big] = \beta (u_x^2 + u_y^2).
\end{align*}

\subsection{Mathematical Difficulties}
There exists a fairly satisfactory well posedness theory for the one dimensional 
equation \eqref{eq:main}.
However, despite its apparent simplicity, the mathematical analysis of \eqref{eq:main} 
is complicated. Independently of the smoothness of the initial data, due to the
nonlinear nature of the equation, singularities may form in the solution
\cite{glassey,ghz1996,ghz1997}.
Therefore, solutions of \eqref{eq:main} should be interpreted in the weak sense:
\begin{definition}
  Set $\Pt=\R \times (0,T)$.  A function
  \begin{equation*}
    u(t,x) \in L^{\infty}\left([0,T];W^{1,p}(\R)\right) \cap C(\Pt), u_t
    \in L^{\infty}\left([0,T];L^{p}(\R)\right),
  \end{equation*}
  for all $p \in [1, 3+q]$, where $q$ is some positive constant, is
  a weak solution of the initial value problem \eqref{eq:main} if it
  satisfies:
  \begin{Definitions}
  \item For all test functions $\test\in \mathcal{D}(\R \times
    [0,T))$ \label{def:w3}
    \begin{equation}
      \label{eq:weaksol}
      \iint_{\Pt} \left( u_t \test_t -c^2(u) u_x \test_x - c(u) c'(u)
        (u_x)^2 \test \right)\,dx \,dt = 0. 
    \end{equation}
  \item $u(\cdot, t) \rightarrow u_0$ in $C\left( [0,T];L^2(\R)
    \right)$ as $t \rightarrow 0^{+}$.%\label{def:w1}
  \item $u_t(\cdot, t) \rightarrow u_1$ as a distribution in $\Pt$
    when $t \rightarrow 0^{+}$.%\label{def:w2}
  \end{Definitions}
\end{definition}
In recent years, there has been an increased interest to understand the
different classes of weak solutions (conservative and dissipative) of the Cauchy
problem \eqref{eq:main}, under the restrictive assumption on the wave speed $c$
(positivity of the derivative of $c$). The literature herein is substantial,
and we will here only give a non-exhaustive overview. Within the existing
framework, we mention the papers by Zhang and Zheng
\cite{zhang1,zhang2,zhang3,zhang4,zhang5,zhang6}, Bressan and Zheng
\cite{bressan} and Holden and Raynaud~\cite{holden}. In fact, taking advantage
of Young measure theory, existence of a global weak solution with initial data
$u_0 \in H^1(\R)$ and $u_1 \in L^2(\R)$ has been proved in \cite{zhang5}.
However, the regularity assumptions on the wave speed $c(u)$ ($c(u)$ is
smooth, bounded, positive with derivative that is non-negative and strictly
positive on the initial data $u_0$) in the analysis of
\cite{zhang1,zhang2,zhang3,zhang4,zhang5,zhang6} precludes consideration of
the physical wave speed given by \eqref{eq:wavespeed}.

A novel approach to the study of \eqref{eq:main} was taken by
Bressan and Zheng \cite{bressan}. 
They have constructed the solutions by introducing new variables
related to the characteristics, leading to a characterization of
singularities in the energy density. The solution $u$, constructed by the above principle,  
is locally Lipschitz continuous
and the map $t \rightarrow u(t,\cdot)$ is continuously differentiable
with values in $L^p_{\mathrm{loc}}(\R)$ for $1 \le p < 2$.

Drawing preliminary motivation from \cite{bressan}, Holden and Raynaud
\cite{holden} provides a rigorous construction of a semigroup of
conservative solutions of \eqref{eq:main}.  Since their construction is based
on energy measures as independent variables, the formation of singularities is
somewhat natural and they were able to overcome the non-physical condition on
wave speed ($c'(u)>0$). Moreover, their analysis can incorporate initial data
$u_0, u_1$ that contain measures.

On the other side, the existence of solutions to two dimensional planar waves
\eqref{eq:main_2d} is completely open. Contrary to its one dimensional
counterpart, it is not possible to rewrite \eqref{eq:main_2d} as a system of
equations in terms of Riemann invariants (for a brief justification, see Sec
~\ref{sec:2d}). Therefore, the same proofs do not apply mutatis mutandis in
the two dimensional case. Having said this, one can of course rewrite
\eqref{eq:main_2d} as a first order system using different change of variables
(see Sec ~\ref{sec:2d}).  However, due to lack of ``symmetry'' of this
formulation, it is hard to establish well posedness of such equations using
this approach. The convergence of numerical schemes (DG or others) to weak
solutions of the 2D equation is also a delicate issue, due to the nonlinearity
associated with the elastic energy. However, in the non-physical one-constant
approximation ($\alpha = \beta$) the equation becomes linear and classical 
convergence results can be applied.

\subsection{Numerical Schemes}
Except under very simplifying assumptions, there does not exist elementary and 
explicit solutions for \eqref{eq:main}.
Moreover, the existence of two classes of weak solutions renders the initial value problem
ill-posed after the formation of singularities. Consequently, robust
numerical schemes are important in the study of the variational wave equation.
Furthermore, capturing conservative solutions
numerically is indeed a delicate issue since we expect that traditional finite difference 
schemes will not yield conservative solutions, due to the intrinsic numerical 
diffusion in these schemes. 

There is a sparsity of efficient numerical schemes for the 1D equation
\eqref{eq:main} available in the literature. We can refer to \cite{ghz1997},
where the authors present some numerical examples to illustrate their theory.
By the way of the theory of Young's measure-valued solutions, Holden et. al.
\cite{hkr2009} proved convergence of the numerical approximation generated by
a semi-discrete finite difference scheme for one-dimensional equation
\eqref{eq:main} to the \emph{dissipative} weak solution of \eqref{eq:main}, under
a restrictive assumption on the wave speed ($c'(u)>0$). 
To overcome such non-physical assumptions, Holden and Raynaud \cite{holden} used 
their analytical construction, as mentioned earlier, to define a numerical method
that can approximate the \emph{conservative} solution. However, the main drawback of this 
method is that it is computationally very expensive as there is no time marching.

Finally, we mention recent papers \cite{koley,Aursand2014Preprint} which deals
with finite difference schemes and discontinuous Galerkin schemes,
respectively, for \eqref{eq:main}. Their main idea was to rewrite
\eqref{eq:main} in the form of a first order systems and design numerical
schemes for those systems.  The key design principle was either energy
conservation or energy dissipation.  In that context, they have presented
schemes that either conserve or dissipate the discrete energy. They also
validated the properties of the schemes via extensive numerical experiments. 
 
Numerical results for the two-dimensional variational wave equation
\eqref{eq:main_2d} are even more sparse than for the one-dimensional case. In
fact, to the best of the authors' knowledge, the only available numerical
experiments are given in the final section of the recent paper by Koley et
al.~\cite{koley}. 

\subsection{Scope and Outline of the Paper}
The purpose of this paper is to develop efficient high-order schemes
for the two-dimensional nonlinear variational wave equation \eqref{eq:main_2d}. 
By using the Discontinuous Galerkin framework we aim to derive
schemes that either \emph{conserve} or \emph{dissipate} a discrete
version of the energy inherited from the variational formulation of the
problem. The proposed DG formulation is in space, and we use high-order
Runge--Kutta schemes to integrate in the temporal dimension.  
Since the behavior of solutions to the 2D equation \eqref{eq:main_2d} is largely
unknown, these schemes will allow us to begin investigate if crucial properties 
of the 1D equation \eqref{eq:main} carry over in the two-dimensional
case. To the best of our knowledge, this is the first systematic numerical
study of the two-dimensional variational wave equation \eqref{eq:main_2d}. 

Our approach for constructing high-order schemes is the RK-DG method
\cite{hill,cockburnlinshu}, where the test and trial functions are
discontinuous piecewise polynomials. In contrast to high order finite-volume
schemes, the high order of accuracy is already built into the finite
dimensional spaces and no reconstruction is needed. Exact or approximate
Riemann solvers from finite volume methods are used to compute the numerical
fluxes between elements. For an energy dissipative scheme we will employ a
combination of dissipative fluxes and, in order to control possible spurious
oscillations near shocks, shock capturing operators
\cite{Johnson1990,chavent,barth}. These methods have recently 
been shown to be entropy
stable for conservation laws \cite{hiltebrand}.
In contrast to for finite volume methods, entropy stability has
gained more attention in finite element methods since one advantage of this
method is that the formulation immediately allows the use of general
unstructured grids.

The shock capturing DG schemes in this paper have the following properties:
\begin{enumerate}
\item The schemes are arbitrarily high-order accurate.
\item The schemes are robust and resolved the solution (including
  possible singularities in the angle $u$) in a stable manner.
\item The energy conservative scheme preserves the discrete energy at the
  semi-discrete level. Using a high-order time stepping method, this property
  also holds in the fully discrete case for all orders of accuracy tested.
\item The energy dissipative scheme dissipates the discrete energy at the
  semi-discrete level. Using a high-order time stepping method, this property
  also holds in the fully discrete case for all orders of accuracy tested.
\end{enumerate}

In the current presentation we consider, for simplicity, a Cartesian grid. The
schemes can however be generalized to more general geometries. For such
applications, it might be useful to write \eqref{eq:main_2d} in the form
\begin{equation}
  u_{tt} - (T(u) \nabla)\left(T(u) \nabla u\right) = 0
  \label{eq:main_2d_alt}
\end{equation}
where
\begin{equation*}
  T(u) = \begin{pmatrix}
    \sqrt{\alpha} \cos(u) & \sqrt{\alpha} \sin(u) \\
    -\sqrt{\beta} \sin(u) & \sqrt{\beta} \cos(u)
  \end{pmatrix}.
\end{equation*}

The rest of the paper is organized as follows: In Section
\ref{sec:2d}, we present energy conservative and energy dissipative
schemes for the one-dimensional equation \eqref{eq:main_2d}. Section
\ref{sec:ham_scheme} 
concerns a first-order Hamiltonian (energy preserving) scheme for comparison.
Section \ref{sec:numerical} contains numerical experiments verifying the order
of convergence, energy stability and efficiency of the schemes. 

\section{Discontinuous Galerkin Schemes in Two-space Dimensions}
\label{sec:2d}
Drawing primary motivation from the one-dimensional case \cite{Aursand2014Preprint}, we aim to design energy conservative and
energy dissipative discontinuous Galerkin schemes of the
two-dimensional version of the nonlinear variational wave equation
\eqref{eq:main_2d}, by rewriting it as a first-order system. First, we briefly mention
why formulation based on Riemann invariants does not work in two dimensional case.

\subsection{The system of equations}
We introduce three new independent variables:
\begin{align*}
  p & := u_t, \\
  v & := \cos(u) u_x + \sin(u) u_y,\\
  w & := \sin(u) u_x - \cos(u) u_y.
\end{align*}
Then, for smooth solutions, we see that
\begin{align*}
v_t &= \cos(u) u_{xt} -\sin(u) u_t u_x + \sin(u) u_{yt} + \cos(u) u_t u_y \\
&=  (\cos(u) u_{t})_x - u_t (\cos(u))_x + (\sin(u) u_{t})_y - u_t (\sin(u))_y -u_t \left(  \sin(u) u_x - \cos(u) u_y \right), 
\end{align*}
and
\begin{align*}
w_t &= \sin(u) u_{xt} +\cos(u) u_t u_x - \cos(u) u_{yt} + \sin(u) u_t u_y \\
&=  (\sin(u) u_{t})_x - u_t (\sin(u))_x - (\cos(u) u_{t})_y + u_t (\cos(u))_y +u_t \left(  \cos(u) u_x + \sin(u) u_y \right).
\end{align*}
Moreover, a straightforward calculation using equation \eqref{eq:main_2d} reveals that
\begin{align*}
p_t  - (\alpha - \beta)& \left( \cos(u) \sin(u) u_x^2 - \cos^2(u) u_x u_y + \sin^2(u) u_x u_y -\cos(u) \sin(u) u_y^2 \right) \\
&= \alpha \left( \cos(u) ( \cos(u) u_x + \sin(u) u_y) \right)_x + \alpha \left( \sin(u) ( \cos(u) u_x + \sin(u) u_y) \right)_y \\
& \qquad + \beta \left( \sin(u) (\sin(u) u_x - \cos(u) u_y ) \right)_x - \beta \left( \cos(u) (\sin(u) u_x - \cos(u) u_y) \right)_y.
\end{align*}
Hence, for smooth solutions, equation \eqref{eq:main_2d} is equivalent
to the following system for $(p,v,w,u)$,
\begin{equation}
  \begin{cases}
    p_t - \alpha (f(u) v)_x -\alpha (g(u) v)_y - \beta (g(u) w)_x + 
    \beta (f(u) w)_y - \alpha v w + \beta v w = 0, &\\
    v_t - (f(u) p)_x + p f(u)_x - (g(u) p)_y + p g(u)_y + pw = 0, &\\
    w_t - (g(u) p)_x + p g(u)_x + (f(u) p)_y - p f(u)_y - pv = 0, & \\
    u_t =p, &
  \end{cases}
  \label{eq:main1_2d}
\end{equation}
where $f(u):= \cos(u)$, and $g(u):= \sin(u)$.
Furthermore, the corresponding energy associated with the system \eqref{eq:main1_2d}
is
\begin{equation}
  \label{eq:energy2_2d}
  \mathcal{E}(t) = \iint_{\R^2} \left( p^2 + \alpha \,v^2 + \beta \, w^2 \right) \,\mathrm{d} x \, \mathrm{d} y.
\end{equation}
A simple calculation shows that smooth solutions of \eqref{eq:main1_2d}
satisfy the energy identity:
\begin{equation}
  \label{eq:egypreserve2_2d}
 \left( p^2 + \alpha \,v^2 + \beta \, w^2 \right)_t + 2 
    \left( \alpha \, p\, f(u) \, v + \beta \, p\, g(u)\, w \right)_x 
  + 2 \left( \alpha \, p\, g(u) \, v - \beta \, p\, f(u)\, w \right)_y=0.
\end{equation}
Hence, the fact that the total energy \eqref{eq:energy2_2d} is conserved
follows from integrating the above identity in space and assuming that
the functions $p,u,v$ and $w$ decay at infinity.

\subsection{The grid}
We begin by introducing some notation needed to define the DG schemes.  Let
the domain $\Omega \subset \mathbb{R}^2$ be decomposed as $\Omega = \cup_{i,j}
\Omega_{ij}$ with $ \Omega_{ij}:= \Omega_i \times \Omega_j$ where $\Omega_i = [x_{i-1/2},x_{i+1/2}]$ and 
$\Omega_j = [y_{j-1/2},y_{j+1/2}]$ for $i,j=1,\cdots,N$. Moreover, we
denote $\Delta x_i = x_{i+1/2} - x_{i-1/2}$ and $ \Delta y_j = y_{j+1/2} - y_{j-1/2}$. Furthermore, we
also denote $x_i = (x_{i-1/2} + x_{i+1/2})/2$ and $ y_j = (y_{j-1/2} + y_{j+1/2})/2$.

Let $u$ be a grid function and denote $u^+_{i+1/2}(y)$ as the
function evaluated at the right side of the cell interface at $x_{i+1/2}$
and let $u^-_{i+1/2}(y)$ denote the
value at the left side. Similarly, we let $u^+_{j+1/2}(x)$ be the
function evaluated at the upper side of the cell interface at $y_{i+1/2}$ 
and let $u^-_{j+1/2}(x)$ denote the value at the lower side. We can then introduce
the jump and, respectively, the average of any grid function $u$ across the
interfaces as
\begin{align*}
  \overline{u}_{i+1/2}(y) &:= \frac{u^+_{i+1/2}(y) + u^-_{i+1/2}(y)}{2},  
  \quad    \overline{u}_{j+1/2}(x) := \frac{u^+_{j+1/2}(x) + u^-_{j+1/2}(x)}{2}, \\
  {\llbracket u \rrbracket}_{i + 1/2}(y) & := u^+_{i+1/2}(y) -u^-_{i+1/2}(y), 
  \quad {\llbracket u \rrbracket}_{j + 1/2}(x)  := u^+_{j+1/2}(x)
  -u^-_{j+1/2}(x).
\end{align*}
Moreover, let $v$ be another grid function. Then the following identities are
readily verified:
\begin{equation}
\begin{aligned}
 \label{eq:useful_2d}
 {\llbracket u v \rrbracket}_{i + 1/2}  = \overline{u}_{i+1/2} \vv_{i +1/2} + { \llbracket u  \rrbracket}_{i + 1/2} \overline{v}_{i+1/2},
\quad {\llbracket u v \rrbracket}_{j + 1/2}  = \overline{u}_{j+1/2} \vv_{j +1/2} + { \llbracket u  \rrbracket}_{j + 1/2} \overline{v}_{j+1/2}
\end{aligned}
\end{equation}

\subsection{Variational Formulation}
We seek an approximation $(p,v, w, u)$ of \eqref{eq:main1_2d} such that for each $t \in [0,T]$,
$p$, $v$, $w$, and $u$ belong to finite dimensional space
\begin{equation*}
  {X}_{h}^s(\Omega) = \left\{ 
    u \in L^2(\Omega) : u|_{\Omega_{ij}} \, \, 
    \text{polynomial of degree} \leq p \right\}.
\end{equation*}

The variational form is derived by multiplying the strong form \eqref{eq:main1_2d}
with test functions $\phi, \nu, \psi, \zeta \in X_{h}^s(\Omega)$ and integrating over 
each element separately. After using integration-by-parts, we obtain
\begin{multline}
  \sum_{i,j=1}^N \int_{\Omega_{ij}} p_t \, \phi \,\mathrm{d} x \, \mathrm{d} y
  + \alpha \sum_{i,j=1}^N \int_{\Omega_{ij}} f(u)\, v \, \phi_x \, \mathrm{d} x \, \mathrm{d} y
  - \alpha \sum_{i,j=1}^N \int_{\Omega_{j}} (f v)_{i+1/2}\,  \phi_{i+1/2}^- \,\mathrm{d} y \\
  + \alpha \sum_{i,j=1}^N \int_{\Omega_{j}} (f v)_{i-1/2} \, \phi_{i-1/2}^+ \,\mathrm{d} y 
  + \alpha \sum_{i,j=1}^N \int_{\Omega_{ij}} g(u)\, v \, \phi_y \, \mathrm{d} x \, \mathrm{d} y
  - \alpha \sum_{i,j=1}^N \int_{\Omega_{i}} (g v)_{j+1/2} \, \phi_{j+1/2}^- \, \mathrm{d} x \\
  + \alpha \sum_{i,j=1}^N \int_{\Omega_{i}} (g v)_{j-1/2} \, \phi_{j-1/2}^+ \, \mathrm{d} x
  + \beta \sum_{i,j=1}^N \int_{\Omega_{ij}} g(u)\, w \, \phi_x \, \mathrm{d} x \, \mathrm{d} y
  - \beta  \sum_{i,j=1}^N \int_{\Omega_{j}} (g w)_{i+1/2} \, \phi_{i+1/2}^- \, \mathrm{d} y\\
  + \beta \sum_{i,j=1}^N \int_{\Omega_{j}} (g w)_{i-1/2} \, \phi_{i-1/2}^+ \, \mathrm{d} y
  - \beta \sum_{i,j=1}^N \int_{\Omega_j} f(u)\, w \, \phi_y \, \mathrm{d} x \, \mathrm{d} y
  + \beta \sum_{i,j=1}^N \int_{\Omega_{i}} (f w)_{j+1/2} \, \phi_{j+1/2}^-  \, \mathrm{d} x \\
  - \beta \sum_{i,j=1}^N \int_{\Omega_{i}} (f w)_{j-1/2} \, \phi_{j-1/2}^+ \, \mathrm{d} x
  - \alpha \sum_{i,j=1}^N \int_{\Omega_{ij}} v\, w \, \phi \, \mathrm{d} x \, \mathrm{d} y
  + \beta \sum_{i,j=1}^N \int_{\Omega_{ij}} v\, w \, \phi \, \mathrm{d} x \, \mathrm{d} y=0, 
  \label{eq:weak_form_vw_a}
\end{multline}
and
\begin{multline}
  \sum_{i,j=1}^N  \int_{\Omega_{ij}} v_t \, \nu \,\mathrm{d} x \, \mathrm{d} y
  + \sum_{i,j=1}^N \int_{\Omega_{ij}} f(u)\, p\, \nu_x \, \mathrm{d} x  \, \mathrm{d} y
  - \sum_{i,j=1}^N \int_{\Omega_j} (f \, p)_{i+1/2} \, \nu_{i+1/2}^- \, \mathrm{d} y \\
  + \sum_{i,j=1}^N \int_{\Omega_j} (f \, p)_{i-1/2} \, \nu_{i-1/2}^+ \, \mathrm{d} y
  - \sum_{i,j=1}^N  \int_{\Omega_{ij}} f(u) \left(p \, \nu\right)_x \, \mathrm{d} x \, \mathrm{d} y
  + \sum_{i,j=1}^N  \int_{\Omega_j} (f)_{i+1/2} \, p_{i+1/2}^- \, \nu_{i+1/2}^- \, \mathrm{d} y\\
  - \sum_{i,j=1}^N \int_{\Omega_j} (f)_{i-1/2} \, p_{i-1/2}^+ \, \nu_{i-1/2}^+  \, \mathrm{d} y
  + \sum_{i,j=1}^N  \int_{\Omega_{ij}} g(u)\, p\, \nu_y \, \mathrm{d} x \, \mathrm{d} y 
  - \sum_{i,j=1}^N \int_{\Omega_{i}} (g \, p)_{j+1/2} \, \nu_{j+1/2}^- \, \mathrm{d} x \\
  + \sum_{i,j=1}^N \int_{\Omega_{i}} (g \, p)_{j-1/2} \, \nu_{j-1/2}^+ \, \mathrm{d} x
  - \sum_{i,j=1}^N  \int_{\Omega_{ij}} g(u) \left(p \, \nu\right)_y \,  \mathrm{d} x \, \mathrm{d} y
  + \sum_{i,j=1}^N  \int_{\Omega_{i}} (g)_{j+1/2} \, p_{j+1/2}^- \, \nu_{j+1/2}^- \, \mathrm{d} x \\
  - \sum_{i,j=1}^N  \int_{\Omega_{i}} (g)_{j-1/2} \, p_{j-1/2}^+ \, \nu_{j-1/2}^+  \, \mathrm{d} x
  + \sum_{i,j=1}^N  \int_{\Omega_{ij}} p\, w \, \nu \, \mathrm{d} x \, \mathrm{d} y=0,
  \label{eq:weak_form_vw_b}
\end{multline}
and
\begin{multline}
  \sum_{i,j=1}^N  \int_{\Omega_{ij}} w_t \, \psi \, \mathrm{d} x \, \mathrm{d} y
  + \sum_{i,j=1}^N  \int_{\Omega_{ij}} g(u)\, p\, \psi_x \, \mathrm{d} x \, \mathrm{d} y
  - \sum_{i,j=1}^N  \int_{\Omega_j} (g \,p)_{i+1/2} \, \psi_{i+1/2}^- \, \mathrm{d} y\\
  + \sum_{i,j=1}^N \int_{\Omega_j} (g \, p)_{i-1/2} \, \psi_{i-1/2}^+  \, \mathrm{d} y
  - \sum_{i,j=1}^N  \int_{\Omega_{ij}} g(u) \left(p \, \psi \right)_x \, \mathrm{d} x \, \mathrm{d} y
  + \sum_{i,j=1}^N  \int_{\Omega_j} (g)_{i+1/2} \, p_{i+1/2}^- \, \psi_{i+1/2}^- \, \mathrm{d} y\\
  - \sum_{i,j=1}^N \int_{\Omega_j} (g)_{i-1/2} \, p_{i-1/2}^+ \, \psi_{i-1/2}^+  \, \mathrm{d} y
  - \sum_{i,j=1}^N  \int_{\Omega_{ij}} f(u)\, p\, \psi_y \, \mathrm{d} x \, \mathrm{d} y
  + \sum_{i,j=1}^N  \int_{\Omega_{i}} (f \, p)_{j+1/2} \, \psi_{j+1/2}^- \, \mathrm{d} x \\
  - \sum_{i,j=1}^N  \int_{\Omega_{i}} (f \, p)_{j-1/2} \, \psi_{j-1/2}^+ \, \mathrm{d} x 
  + \sum_{i,j=1}^N  \int_{\Omega_{ij}} f(u) \left(p \, \psi\right)_y \mathrm{d} x \, \mathrm{d} y
  - \sum_{i,j=1}^N   \int_{\Omega_{i}} (f)_{j+1/2} \, p_{j+1/2}^- \, \psi_{j+1/2}^- \, \mathrm{d} x  \\
  + \sum_{i,j=1}^N  \int_{\Omega_{i}}  (f)_{j-1/2} \, p_{j-1/2}^+ \, \psi_{j-1/2}^+  \, \mathrm{d} x 
  - \sum_{i,j=1}^N  \int_{\Omega_{ij}} p\, v \, \psi \, \mathrm{d} x \, \mathrm{d} y =0,
  \label{eq:weak_form_vw_c}
\end{multline}
and 
\begin{align}
 \sum_{i,j=1}^N  \int_{\Omega_{ij}} u_t \, \zeta \, \mathrm{d} x \, \mathrm{d} y = \sum_{i,j=1}^N
  \int_{\Omega_{ij}} p \, \zeta \, \mathrm{d} x \, \mathrm{d} y.
  \label{eq:weak_form_vw_d}
\end{align}

\begin{remark}
  Admittedly, the notation used in
  \eqref{eq:weak_form_vw_a}--\eqref{eq:weak_form_vw_d} is more cumbersome than
  the vector notation often seen in the DG literature. The purpose of this is
  to be able to treat the fluxes in the different equations differently in
  order to ensure energy conservation. Also, since the proposed scheme is for
  the nonlinear variational wave equation, not a general class of 
  wave equations, we hope to avoid unnecessary confusion by writing fluxes explicitly.
\end{remark}

In order to complete the description of the above schemes, we need to specify numerical flux functions. 

\subsection{Energy Preserving Scheme}

For a conservative scheme, we use the central numerical flux
\begin{equation*}
  (f)_{k \pm 1/2} = \overline{f}_{k \pm 1/2} \quad \text{and} \quad (f g)_{k \pm 1/2} 
      = \overline{f}_{k \pm 1/2} \overline{g}_{k \pm 1/2},
\end{equation*}
for any grid functions $f,g \in X^s_h(\Omega)$.
An energy preserving (spatial) DG scheme based on the weak formulation
\eqref{eq:weak_form_vw_a}--\eqref{eq:weak_form_vw_d} becomes: 
Find $p, v, w, u \in X_{h}^s(\Omega)$ such that
\begin{multline}
  \sum_{i,j=1}^N \int_{\Omega_{ij}} p_t \, \phi \,\mathrm{d} x \, \mathrm{d} y
  + \alpha \sum_{i,j=1}^N \int_{\Omega_{ij}} f(u)\, v \, \phi_x \, \mathrm{d} x \, \mathrm{d} y
  - \alpha \sum_{i,j=1}^N \int_{\Omega_{j}} \overline{f}_{i+1/2} \, \overline{v}_{i+1/2} \,  \phi_{i+1/2}^- \,\mathrm{d} y \\
  + \alpha \sum_{i,j=1}^N \int_{\Omega_{j}} \overline{f}_{i-1/2} \, \overline{v}_{i-1/2} \, \phi_{i-1/2}^+ \,\mathrm{d} y 
  + \alpha \sum_{i,j=1}^N \int_{\Omega_{ij}} g(u)\, v \, \phi_y \, \mathrm{d} x \, \mathrm{d} y
  - \alpha \sum_{i,j=1}^N \int_{\Omega_{i}} \overline{g}_{j+1/2} \,\overline{v}_{j+1/2} \, \phi_{j+1/2}^- \, \mathrm{d} x \\
  + \alpha \sum_{i,j=1}^N \int_{\Omega_{i}} \overline{g}_{j-1/2} \,\overline{v}_{j-1/2} \, \phi_{j-1/2}^+ \, \mathrm{d} x
  + \beta \sum_{i,j=1}^N \int_{\Omega_{ij}} g(u)\, w \, \phi_x \, \mathrm{d} x \, \mathrm{d} y
  - \beta  \sum_{i,j=1}^N \int_{\Omega_{j}} \overline{g}_{i+1/2} \, \overline{w}_{i+1/2} \, \phi_{i+1/2}^- \, \mathrm{d} y\\
  + \beta \sum_{i,j=1}^N \int_{\Omega_{j}} \overline{g}_{i-1/2} \, \overline{w}_{i-1/2} \, \phi_{i-1/2}^+ \, \mathrm{d} y
  - \beta \sum_{i,j=1}^N \int_{\Omega_j} f(u)\, w \, \phi_y \, \mathrm{d} x \, \mathrm{d} y
  + \beta \sum_{i,j=1}^N \int_{\Omega_{i}}  \overline{f}_{j+1/2} \, \overline{w}_{j+1/2} \, \phi_{j+1/2}^-  \, \mathrm{d} x \\
  - \beta \sum_{i,j=1}^N \int_{\Omega_{i}}  \overline{f}_{j-1/2} \, \overline{w}_{j-1/2} \, \phi_{j-1/2}^+ \, \mathrm{d} x
  - \alpha \sum_{i,j=1}^N \int_{\Omega_{ij}} v\, w \, \phi \, \mathrm{d} x \, \mathrm{d} y
  + \beta \sum_{i,j=1}^N \int_{\Omega_{ij}} v\, w \, \phi \, \mathrm{d} x \, \mathrm{d} y=0, 
  \label{eq:weak_form_vw_a_con}
\end{multline}
for all $\phi \in X_{\Delta x}^s(\Omega)$,
\begin{multline}
  \sum_{i,j=1}^N  \int_{\Omega_{ij}} v_t \, \nu \,\mathrm{d} x \, \mathrm{d} y
  + \sum_{i,j=1}^N \int_{\Omega_{ij}} f(u)\, p\, \nu_x \, \mathrm{d} x  \, \mathrm{d} y
  - \sum_{i,j=1}^N \int_{\Omega_j} \overline{f}_{i+1/2} \, \overline{p}_{i+1/2} \, \nu_{i+1/2}^- \, \mathrm{d} y \\
  + \sum_{i,j=1}^N \int_{\Omega_j} \overline{f}_{i-1/2} \, \overline{p}_{i-1/2} \, \nu_{i-1/2}^+ \, \mathrm{d} y
  - \sum_{i,j=1}^N  \int_{\Omega_{ij}} f(u) \left(p \, \nu\right)_x \, \mathrm{d} x \, \mathrm{d} y
  + \sum_{i,j=1}^N  \int_{\Omega_j} \overline{f}_{i+1/2}\, p_{i+1/2}^- \, \nu_{i+1/2}^- \, \mathrm{d} y\\
  - \sum_{i,j=1}^N \int_{\Omega_j} \overline{f}_{i-1/2} \, p_{i-1/2}^+ \, \nu_{i-1/2}^+  \, \mathrm{d} y
  + \sum_{i,j=1}^N  \int_{\Omega_{ij}} g(u)\, p\, \nu_y \, \mathrm{d} x \, \mathrm{d} y 
  - \sum_{i,j=1}^N \int_{\Omega_{i}} \overline{g}_{j+1/2} \, \overline{p}_{j+1/2} \, \nu_{j+1/2}^- \, \mathrm{d} x \\
  + \sum_{i,j=1}^N \int_{\Omega_{i}} \overline{g}_{j-1/2} \, \overline{p}_{j-1/2} \, \nu_{j-1/2}^+ \, \mathrm{d} x
  - \sum_{i,j=1}^N  \int_{\Omega_{ij}} g(u) \left(p \, \nu\right)_y \,  \mathrm{d} x \, \mathrm{d} y
  + \sum_{i,j=1}^N  \int_{\Omega_{i}}  \overline{g}_{j+1/2} \, p_{j+1/2}^- \, \nu_{j+1/2}^- \, \mathrm{d} x \\
  - \sum_{i,j=1}^N  \int_{\Omega_{i}}  \overline{g}_{j-1/2} \, p_{j-1/2}^+ \, \nu_{j-1/2}^+  \, \mathrm{d} x
  + \sum_{i,j=1}^N  \int_{\Omega_{ij}} p\, w \, \nu \, \mathrm{d} x \, \mathrm{d} y=0,
  \label{eq:weak_form_vw_b_con}
\end{multline}
for all $\nu \in X_{h}^s(\Omega)$,
\begin{multline}
  \sum_{i,j=1}^N  \int_{\Omega_{ij}} w_t \, \psi \, \mathrm{d} x \, \mathrm{d} y
  + \sum_{i,j=1}^N  \int_{\Omega_{ij}} g(u)\, p\, \psi_x \, \mathrm{d} x \, \mathrm{d} y
  - \sum_{i,j=1}^N  \int_{\Omega_j} \overline{g}_{i+1/2} \, \overline{p}_{i+1/2}  \, \psi_{i+1/2}^- \, \mathrm{d} y\\
  + \sum_{i,j=1}^N \int_{\Omega_j} \overline{g}_{i-1/2} \, \overline{p}_{i-1/2}  \, \psi_{i-1/2}^+  \, \mathrm{d} y
  - \sum_{i,j=1}^N  \int_{\Omega_{ij}} g(u) \left(p \, \psi \right)_x \, \mathrm{d} x \, \mathrm{d} y
  + \sum_{i,j=1}^N  \int_{\Omega_j} \overline{g}_{i+1/2} \, p_{i+1/2}^- \, \psi_{i+1/2}^- \, \mathrm{d} y\\
  - \sum_{i,j=1}^N \int_{\Omega_j} \overline{g}_{i-1/2} \, p_{i-1/2}^+ \, \psi_{i-1/2}^+  \, \mathrm{d} y
  - \sum_{i,j=1}^N  \int_{\Omega_{ij}} f(u)\, p\, \psi_y \, \mathrm{d} x \, \mathrm{d} y
  + \sum_{i,j=1}^N  \int_{\Omega_{i}} \overline{f}_{j+1/2} \, \overline{p}_{j+1/2} \, \psi_{j+1/2}^- \, \mathrm{d} x \\
  - \sum_{i,j=1}^N  \int_{\Omega_{i}} \overline{f}_{j-1/2} \, \overline{p}_{j-1/2} \, \psi_{j-1/2}^+ \, \mathrm{d} x 
  + \sum_{i,j=1}^N  \int_{\Omega_{ij}} f(u) \left(p \, \psi\right)_y \mathrm{d} x \, \mathrm{d} y
  - \sum_{i,j=1}^N   \int_{\Omega_{i}} \overline{f}_{j+1/2} \, p_{j+1/2}^- \, \psi_{j+1/2}^- \, \mathrm{d} x  \\
  + \sum_{i,j=1}^N  \int_{\Omega_{i}}  \overline{f}_{j-1/2} \, p_{j-1/2}^+ \, \psi_{j-1/2}^+  \, \mathrm{d} x 
  - \sum_{i,j=1}^N  \int_{\Omega_{ij}} p\, v \, \psi \, \mathrm{d} x \, \mathrm{d} y =0,
  \label{eq:weak_form_vw_c_con}
\end{multline}
for all $\psi \in X_{h}^s(\Omega)$ and 
\begin{align}
 \sum_{i,j=1}^N  \int_{\Omega_{ij}} u_t \, \zeta \, \mathrm{d} x \, \mathrm{d} y = \sum_{i,j=1}^N
  \int_{\Omega_{ij}} p \, \zeta \, \mathrm{d} x \, \mathrm{d} y.
  \label{eq:weak_form_vw_d_con}
\end{align}
for all $\zeta \in X_{h}^s(\Omega)$.

The above scheme preserves a discrete version of the energy, as shown
in the following theorem:
\begin{theorem} \label{thm:energy_preservation}
  Let $p$, $v$ and $w$ be approximate
  solutions generated by the scheme 
  \eqref{eq:weak_form_vw_a_con}--\eqref{eq:weak_form_vw_d_con} with periodic
  boundary conditions. Then
  \begin{align*}
    \frac{d}{dt} \sum_{i,j=1}^N \int_{\Omega_{ij}} \left(p^2(t)
      + \alpha \, v^2(t) + \beta \, w^2(t)) \right) \, \mathrm{d} x \, \mathrm{d} y=0.
  \end{align*}
\end{theorem}
  
\begin{proof}
Let $p$, $v$ and $w$ be numerical solutions generated by the scheme
\eqref{eq:weak_form_vw_a_con}--\eqref{eq:weak_form_vw_d_con}. Since those
equations hold for any $\phi, \nu, \psi \in X_{h}^s(\Omega)$, they hold
in particular for $\phi =p,  \nu= v $ and $\psi=w$. We can then calculate
\begin{align*}
\frac{d}{dt} &  \sum_{i,j=1}^N  \int_{\Omega_{ij}} \left(p^2(t) + \alpha \, v^2(t) + \beta \, w^2(t)) \right) \,\mathrm{d} x \, \mathrm{d} y
=2 \sum_{i,j=1}^N \int_{\Omega_{ij}} \left(pp_t + \alpha \, vv_t + \beta \, ww_t \right) \,\mathrm{d} x \, \mathrm{d} y \\
&= 2 \alpha \sum_{i,j=1}^N \int_{\Omega_{j}}  \overline{f}_{i+1/2} \left( \overline{v}_{i+1/2} \jump{p}_{i+1/2} 
     + \overline{p}_{i+1/2} \jump{v}_{i+1/2} - \jump{p v}_{i+1/2}\right) \,\mathrm{d} y \\
& +2 \alpha \sum_{i,j=1}^N \int_{\Omega_{i}}  \overline{g}_{j+1/2} \left( \overline{v}_{j+1/2} \jump{p}_{j+1/2} 
     + \overline{p}_{j+1/2} \jump{v}_{j+1/2} - \jump{p v}_{j+1/2}\right)\,\mathrm{d} x \\
& +2 \beta \sum_{i,j=1}^N \int_{\Omega_{j}}  \overline{g}_{i+1/2} \left( \overline{v}_{i+1/2} \jump{p}_{i+1/2} 
     + \overline{p}_{i+1/2} \jump{w}_{i+1/2} - \jump{p w}_{i+1/2}\right) \,\mathrm{d} y \\
& +2 \alpha \sum_{i,j=1}^N \int_{\Omega_{i}}  \overline{f}_{j+1/2} \left( - \overline{w}_{j+1/2} \jump{p}_{j+1/2} 
     - \overline{p}_{j+1/2} \jump{w}_{j+1/2} + \jump{p w}_{j+1/2}\right)\,\mathrm{d} x =0,
\end{align*}
where we have used the periodic boundary conditions and the 
   identities \eqref{eq:useful_2d}.
\end{proof}

\begin{remark} 
  Theorem \ref{thm:energy_preservation} and similar results to follow
  explicitly assume periodic boundary conditions. It is however
  straightforward to show that these results also hold for certain other
  situations such as with compactly supported or decaying data.
\end{remark}

\subsection{Energy Dissipating Scheme}
Note that the above designed energy conservative scheme
\eqref{eq:weak_form_vw_a_con}--\eqref{eq:weak_form_vw_d_con} is expected to 
approximate a conservative solution of the underlying system \eqref{eq:main_2d}. 
To attempt to approximate a dissipative solution of \eqref{eq:main_2d}, one has to 
add \emph{numerical viscosity}. In this work we propose adding viscosity in
the numerical fluxes (scaled by the maximum wave speed) as well as a 
\emph{shock capturing operator} dissipating energy near shocks or discontinuities. 
Specifically, we propose the following modification of the energy conservative
scheme \eqref{eq:weak_form_vw_a_con}--\eqref{eq:weak_form_vw_d_con}:

Denoting
\begin{equation*}
  s_{i \pm1/2} = \max\{c^-_{i \pm1/2},c^+_{i \pm 1/2}\} \, \, \text{and} 
  \,\, s_{j \pm 1/2} = \max\{b^-_{j \pm 1/2},b^+_{j \pm 1/2}\}
\end{equation*} 
for the maximal local wave velocity,
a dissipative version of the DG scheme is then given by the following:
Find $p, v, w, u \in X_{h}^s(\Omega)$ 
such that
\begin{multline}
  \sum_{i,j=1}^N \int_{\Omega_{ij}} p_t \, \phi \,\mathrm{d} x \, \mathrm{d} y
  + \alpha \sum_{i,j=1}^N \int_{\Omega_{ij}} f(u)\, v \, \phi_x \, \mathrm{d} x \, \mathrm{d} y \\
  \qquad  -\alpha \sum_{i,j=1}^N \int_{\Omega_{j}} 
  \underbrace{\left(\overline{f}_{i+1/2} \, \overline{v}_{i+1/2} + \frac{1}{2}  s_{i+1/2}
  \jump{p}_{i+1/2}  \right)}_{\text{diffusive flux in x-direction}}\,  \phi_{i+1/2}^- \,\mathrm{d} y \\
  \qquad \qquad + \alpha \sum_{i,j=1}^N \int_{\Omega_{j}} 
   \underbrace{\left(\overline{f}_{i-1/2} \, \overline{v}_{i-1/2} + \frac{1}{2}  s_{i-1/2} \jump{p}_{i-1/2}  \right)}_{\text{diffusive flux in x-direction}}\, \phi_{i-1/2}^+ \,\mathrm{d} y \\
  + \alpha \sum_{i,j=1}^N \int_{\Omega_{ij}} g(u)\, v \, \phi_y \, \mathrm{d} x \, \mathrm{d} y
  - \alpha \sum_{i,j=1}^N \int_{\Omega_{i}} 
   \underbrace{\left(\overline{g}_{j+1/2} \,\overline{v}_{j+1/2} + \frac{1}{2}  s_{j+1/2} \jump{p}_{j+1/2}   \right)}_{\text{diffusive flux in y-direction}}\, \phi_{j+1/2}^- \, \mathrm{d} x \\
  + \alpha \sum_{i,j=1}^N \int_{\Omega_{i}} 
   \underbrace{\left(\overline{g}_{j-1/2} \,\overline{v}_{j-1/2} + \frac{1}{2}  s_{j-1/2} \jump{p}_{j-1/2}   \right)}_{\text{diffusive flux in y-direction}}\, \phi_{j-1/2}^+ \, \mathrm{d} x 
  + \beta \sum_{i,j=1}^N \int_{\Omega_{ij}} g(u)\, w \, \phi_x \, \mathrm{d} x \, \mathrm{d} y \\
  - \beta  \sum_{i,j=1}^N \int_{\Omega_{j}} \overline{g}_{i+1/2} \, \overline{w}_{i+1/2} \, \phi_{i+1/2}^- \, \mathrm{d} y
  + \beta \sum_{i,j=1}^N \int_{\Omega_{j}} \overline{g}_{i-1/2} \, \overline{w}_{i-1/2} \, \phi_{i-1/2}^+ \, \mathrm{d} y
  - \beta \sum_{i,j=1}^N \int_{\Omega_j} f(u)\, w \, \phi_y \, \mathrm{d} x \, \mathrm{d} y \\
  + \beta \sum_{i,j=1}^N \int_{\Omega_{i}}  \overline{f}_{j+1/2} \, \overline{w}_{j+1/2} \, \phi_{j+1/2}^-  \, \mathrm{d} x 
  - \beta \sum_{i,j=1}^N \int_{\Omega_{i}}  \overline{f}_{j-1/2} \, \overline{w}_{j-1/2} \, \phi_{j-1/2}^+ \, \mathrm{d} x
  - \alpha \sum_{i,j=1}^N \int_{\Omega_{ij}} v\, w \, \phi \, \mathrm{d} x \, \mathrm{d} y \\
  + \beta \sum_{i,j=1}^N \int_{\Omega_{ij}} v\, w \, \phi \, \mathrm{d} x \, \mathrm{d} y
  = - \underbrace{\sum_{i,j=1}^N \varepsilon_{ij} \int_{\Omega_{ij}} \left(p_x \,\phi_x + p_y \,\phi_y \right) \, \mathrm{d} x \, \mathrm{d} y }_{\text{shock capturing operator}},  
  \label{eq:weak_form_vw_a_diss}
\end{multline}
for all $\phi \in X_{h}^s(\Omega)$,
\begin{multline}
  \sum_{i,j=1}^N  \int_{\Omega_{ij}} v_t \, \nu \,\mathrm{d} x \, \mathrm{d} y
  + \sum_{i,j=1}^N \int_{\Omega_{ij}} f(u)\, p\, \nu_x \, \mathrm{d} x  \, \mathrm{d} y
  - \sum_{i,j=1}^N \int_{\Omega_j} 
  \underbrace{\left(\overline{f}_{i+1/2} \, \overline{p}_{i+1/2} + \frac{1}{2}  s_{i+1/2} \jump{v}_{i+1/2} \right)}_{\text{diffusive flux in x-direction}}\, \nu_{i+1/2}^- \, \mathrm{d} y \\
  + \sum_{i,j=1}^N \int_{\Omega_j} 
  \underbrace{\left(\overline{f}_{i-1/2} \, \overline{p}_{i-1/2}  + \frac{1}{2}  s_{i-1/2} \jump{v}_{i-1/2} \right)}_{\text{diffusive flux in x-direction}}\, \nu_{i-1/2}^+ \, \mathrm{d} y \\
  - \sum_{i,j=1}^N  \int_{\Omega_{ij}} f(u) \left(p \, \nu\right)_x \, \mathrm{d} x \, \mathrm{d} y
  + \sum_{i,j=1}^N  \int_{\Omega_j} \overline{f}_{i+1/2}\, p_{i+1/2}^- \, \nu_{i+1/2}^- \, \mathrm{d} y\\
  - \sum_{i,j=1}^N \int_{\Omega_j} \overline{f}_{i-1/2} \, p_{i-1/2}^+ \, \nu_{i-1/2}^+  \, \mathrm{d} y
  + \sum_{i,j=1}^N  \int_{\Omega_{ij}} g(u)\, p\, \nu_y \, \mathrm{d} x \, \mathrm{d} y \\
  - \sum_{i,j=1}^N \int_{\Omega_{i}} 
  \underbrace{\left(\overline{g}_{j+1/2} \, \overline{p}_{j+1/2} + \frac{1}{2} s_{j+1/2} \jump{v}_{j+1/2}  \right)}_{\text{diffusive flux in y-direction}} \, \nu_{j+1/2}^- \, \mathrm{d} x \\
  + \sum_{i,j=1}^N \int_{\Omega_{i}} 
  \underbrace{\left(\overline{g}_{j-1/2} \, \overline{p}_{j-1/2} + \frac{1}{2}  s_{j-1/2} \jump{v}_{j-1/2}\right)}_{\text{diffusive flux in y-direction}}\, \nu_{j-1/2}^+ \, \mathrm{d} x \\
  - \sum_{i,j=1}^N  \int_{\Omega_{ij}} g(u) \left(p \, \nu\right)_y \,  \mathrm{d} x \, \mathrm{d} y
  + \sum_{i,j=1}^N  \int_{\Omega_{i}}  \overline{g}_{j+1/2} \, p_{j+1/2}^- \, \nu_{j+1/2}^- \, \mathrm{d} x \\
  - \sum_{i,j=1}^N  \int_{\Omega_{i}}  \overline{g}_{j-1/2} \, p_{j-1/2}^+ \, \nu_{j-1/2}^+  \, \mathrm{d} x
  + \sum_{i,j=1}^N  \int_{\Omega_{ij}} p\, w \, \nu \, \mathrm{d} x \, \mathrm{d} y
  =- \underbrace{\sum_{i,j=1}^N \varepsilon_{ij} \int_{\Omega_{ij}} \left(v_x \,\nu_x + v_y \,\nu_y \right) \, \mathrm{d} x \, \mathrm{d} y }_{\text{shock capturing operator}},
  \label{eq:weak_form_vw_b_diss}
\end{multline}
for all $\nu \in X_{h}^s(\Omega)$,
\begin{multline}
  \sum_{i,j=1}^N  \int_{\Omega_{ij}} w_t \, \psi \, \mathrm{d} x \, \mathrm{d} y
  + \sum_{i,j=1}^N  \int_{\Omega_{ij}} g(u)\, p\, \psi_x \, \mathrm{d} x \, \mathrm{d} y \\
  - \sum_{i,j=1}^N  \int_{\Omega_j} 
  \underbrace{\left(\overline{g}_{i+1/2} \, \overline{p}_{i+1/2} + \frac{1}{2}  s_{i+1/2} \jump{w}_{i+1/2} \right)}_{\text{diffusive flux in x-direction}} \, \psi_{i+1/2}^- \,\mathrm{d} y\\
  + \sum_{i,j=1}^N \int_{\Omega_j} 
  \underbrace{\left(\overline{g}_{i-1/2} \, \overline{p}_{i-1/2} + \frac{1}{2}  s_{i-1/2} \jump{w}_{i-1/2}\right)}_{\text{diffusive flux in x-direction}} \, \psi_{i-1/2}^+  \, \mathrm{d} y\\
  - \sum_{i,j=1}^N  \int_{\Omega_{ij}} g(u) \left(p \, \psi \right)_x \, \mathrm{d} x \, \mathrm{d} y
  + \sum_{i,j=1}^N  \int_{\Omega_j} \overline{g}_{i+1/2} \, p_{i+1/2}^- \, \psi_{i+1/2}^- \, \mathrm{d} y\\
  - \sum_{i,j=1}^N \int_{\Omega_j} \overline{g}_{i-1/2} \, p_{i-1/2}^+ \, \psi_{i-1/2}^+  \, \mathrm{d} y
  - \sum_{i,j=1}^N  \int_{\Omega_{ij}} f(u)\, p\, \psi_y \, \mathrm{d} x \, \mathrm{d} y \\
  + \sum_{i,j=1}^N  \int_{\Omega_{i}} 
  \underbrace{\left(\overline{f}_{j+1/2} \, \overline{p}_{j+1/2} - \frac{1}{2}  s_{j+1/2} \jump{w}_{j+1/2}\right)}_{\text{diffusive flux in y-direction}} \, \psi_{j+1/2}^- \,\mathrm{d} x \\
  - \sum_{i,j=1}^N  \int_{\Omega_{i}} 
  \underbrace{\left(\overline{f}_{j-1/2} \, \overline{p}_{j-1/2} - \frac{1}{2}  s_{j-1/2} \jump{w}_{j-1/2} \right)}_{\text{diffusive flux in y-direction}}\, \psi_{j-1/2}^+ \, \mathrm{d} x \\
  + \sum_{i,j=1}^N  \int_{\Omega_{ij}} f(u) \left(p \, \psi\right)_y \mathrm{d} x \, \mathrm{d} y
  - \sum_{i,j=1}^N   \int_{\Omega_{i}} \overline{f}_{j+1/2} \, p_{j+1/2}^- \, \psi_{j+1/2}^- \, \mathrm{d} x  \\
  + \sum_{i,j=1}^N  \int_{\Omega_{i}}  \overline{f}_{j-1/2} \, p_{j-1/2}^+ \, \psi_{j-1/2}^+  \, \mathrm{d} x 
  - \sum_{i,j=1}^N  \int_{\Omega_{ij}} p\, v \, \psi \, \mathrm{d} x \, \mathrm{d} y 
  =- \underbrace{\sum_{i,j=1}^N \varepsilon_{ij} \int_{\Omega_{ij}} \left(w_x \,\psi_x + w_y \,\psi_y \right) \, \mathrm{d} x \, \mathrm{d} y }_{\text{shock capturing operator}},
  \label{eq:weak_form_vw_c_diss}
\end{multline}
for all $\psi \in X_{h}^s(\Omega)$, 
\begin{align}
 \sum_{i,j=1}^N  \int_{\Omega_{ij}} u_t \, \zeta \, \mathrm{d} x \, \mathrm{d} y = \sum_{i,j=1}^N
  \int_{\Omega_{ij}} p \, \zeta \, \mathrm{d} x \, \mathrm{d} y.
  \label{eq:weak_form_vw_d_diss}
\end{align}
for all $\zeta \in X_{h}^s(\Omega)$.

The scaling parameter $\varepsilon$ in the \emph{shock capturing operator} is given by
\begin{equation}
  \varepsilon_{ij} = \frac{h_{ij} \, C \, \overline{\mathrm{Res}}}
  {\left(\int_{\Omega_{ij}}(p_x^2 + v_x^2 + w_x^2) \mathrm{d} x \, \mathrm{d} y
  + \int_{\Omega_{ij}}(p_y^2 + v_y^2 + w_y^2) \mathrm{d} x \, \mathrm{d} y 
   \right)^{1/2}+ h_{ij}^\theta  } 
  \label{eq:epsilon_RS}
\end{equation}
where $C > 0$ is a constant, $\theta \geq 1/2$, $h_{ij} = \max\{\Delta x_{ij},\Delta
y_{ij} \}$ and
\begin{equation} \label{eq:residual_rms}
  \overline{\mathrm{Res}} = \left(\int_{\Omega_{ij}} (\mathrm{Res})^2 
  \mathrm{d} x \, \mathrm{d} y \right)^{1/2}
\end{equation}
with
\begin{equation} \label{eq:residual}
  \mathrm{Res} = \left( p^2 + \alpha \,v^2 + \beta \, w^2 \right)_t + \left( \alpha \, p\, f(u) \, v + \beta \, p\, g(u)\, w \right)_x 
  + \left( \alpha \, p\, g(u) \, v - \beta \, p\, f(u)\, w \right)_y.
\end{equation}
The rationale for the scaling parameter is as follows: For smooth solutions of
\eqref{eq:main1_2d} the conservation law \eqref{eq:egypreserve2_2d} is
fulfilled. The numerical solution is then expected to fulfill the same
conservation law up to the spatial and temporal accuracy of the scheme. The shock
capturing operator will therefore vanish in smooth regions, while introducing
added dissipation near shocks and discontinuities.

The above scheme dissipates a discrete version of the energy, as shown in the following 
theorem:
\begin{theorem}
  Let $p$, $v$ and $w$ be approximate
  solutions generated by the scheme 
  \eqref{eq:weak_form_vw_a_con}--\eqref{eq:weak_form_vw_d_con} with periodic
  boundary conditions. Then
  \begin{align*}
    \frac{d}{dt} \sum_{i,j=1}^N \int_{\Omega_{ij}} \left(p^2(t)
      + \alpha \, v^2(t) + \beta \, w^2(t)) \right) \, \mathrm{d} x \,
      \mathrm{d} y \leq 0.
  \end{align*}
\end{theorem}

\begin{proof}
  By using the result from Theorem \ref{thm:energy_preservation}, we can write
  \begin{equation}
  \begin{aligned}
    \frac{d}{dt} \sum_{i,j=1}^N \int_{\Omega_{ij}} &\left(p^2(t)
      + \alpha \, v^2(t) + \beta \, w^2(t)) \right) \, \mathrm{d} x \,
      \mathrm{d} y 
      \\ &= 
      -2 \sum_{i,j=1}^N \eps_{ij} \int_{\Omega_{ij}} \left( p_x^2 + p_y^2
      +v_x^2 + v_y^2 +w_x^2 + w_y^2\right) \mathrm{d} x \, \mathrm{d} y \\
      &+ \alpha \sum_{i,j=1}^N \int_{\Omega_j} \left( s_{i+1/2} \jump{p}_{i+1/2}
      p_{i+1/2}^- - s_{i-1/2} \jump{p}_{i-1/2} p_{i-1/2}^+ \right) \mathrm{d} y
      \\
      &+ \alpha \sum_{i,j=1}^N \int_{\Omega_i} \left( s_{j+1/2} \jump{p}_{j+1/2}
      p_{j+1/2}^- - s_{j-1/2} \jump{p}_{j-1/2} p_{j-1/2}^+ \right) \mathrm{d} x
      \\
      &+ \alpha \sum_{i,j=1}^N \int_{\Omega_j} \left( s_{i+1/2} \jump{v}_{i+1/2}
      v_{i+1/2}^- - s_{i-1/2} \jump{v}_{i-1/2} v_{i-1/2}^+ \right) \mathrm{d} y
      \\
      &+ \alpha \sum_{i,j=1}^N \int_{\Omega_i} \left( s_{j+1/2} \jump{v}_{j+1/2}
      v_{j+1/2}^- - s_{j-1/2} \jump{v}_{j-1/2} v_{j-1/2}^+ \right) \mathrm{d} x
      \\
      &+ \beta \sum_{i,j=1}^N \int_{\Omega_j} \left( s_{i+1/2} \jump{w}_{i+1/2}
      w_{i+1/2}^- - s_{i-1/2} \jump{w}_{i-1/2} w_{i-1/2}^+ \right) \mathrm{d} y
      \\
      &+ \beta \sum_{i,j=1}^N \int_{\Omega_i} \left( s_{j+1/2} \jump{w}_{j+1/2}
      w_{j+1/2}^- - s_{j-1/2} \jump{w}_{j-1/2} w_{j-1/2}^+ \right) \mathrm{d} x
  \end{aligned}
  \end{equation}
  Now, since the periodic boundary condition lends the relation
  \begin{equation}
    \sum_{i,j=1}^N \left(s_{i+1/2} \jump{a}_{i+1/2} a_{i+1/2}^-
                         - s_{i-1/2} \jump{a}_{i-1/2} a_{i-1/2}^+ \right)
                 = - \sum_{i.j=1}^N s_{i+1/2} \jump{a}_{i+1/2}^2,
  \end{equation}
  we can write
  \begin{equation}
  \begin{aligned}
    \frac{d}{dt} \sum_{i,j=1}^N \int_{\Omega_{ij}} &\left(p^2(t)
      + \alpha \, v^2(t) + \beta \, w^2(t)) \right) \, \mathrm{d} x \,
      \mathrm{d} y 
      \\ &= 
      -2 \sum_{i,j=1}^N \eps_{ij} \int_{\Omega_{ij}} \left( p_x^2 + p_y^2
      +v_x^2 + v_y^2 +w_x^2 + w_y^2\right) \mathrm{d} x \, \mathrm{d} y \\
      &- \alpha \sum_{i,j=1}^N \int_{\Omega_j} s_{i+1/2}
      \jump{p}_{i+1/2}^2  \mathrm{d} y
      - \alpha \sum_{i,j=1}^N \int_{\Omega_i} s_{j+1/2} \jump{p}_{j+1/2}^2
      \mathrm{d} x
      \\
      &- \alpha \sum_{i,j=1}^N \int_{\Omega_j} s_{i+1/2} \jump{v}_{i+1/2}^2
      \mathrm{d} y
      - \alpha \sum_{i,j=1}^N \int_{\Omega_i} s_{j+1/2} \jump{v}_{j+1/2}^2
      \mathrm{d} x
      \\
      &- \beta \sum_{i,j=1}^N \int_{\Omega_j} s_{i+1/2} \jump{w}_{i+1/2}^2
      \mathrm{d} y
      - \beta \sum_{i,j=1}^N \int_{\Omega_i} s_{j+1/2} \jump{w}_{j+1/2}^2
      \mathrm{d} x.
  \end{aligned}
  \end{equation}
  The result then follows from the positivity of $\eps_{ij}$, $s$ and the physical
  parameters $\alpha$ and $\beta$.
\end{proof}

\section{Energy Preserving Scheme Based On a Variational Formulation}
\label{sec:ham_scheme}
It is worth noting that all the previous schemes were designed by rewriting
the variational wave equation \eqref{eq:main_2d} as first-order systems and
approximating these systems. However, one can also design a scheme for the
original variational wave equation \eqref{eq:main_2d}. To achieve this, we
design an energy conservative scheme by approximating the nonlinear wave
equation \eqref{eq:main_2d} directly. We proceed by rewriting the nonlinear
wave equation \eqref{eq:main_2d} in the general form:
\begin{equation}
  \label{eq:vari}
  u_{tt} = - \frac{\delta H}{\delta u},
\end{equation}
with 
$$
H = H(u, u_x, u_y):= \frac{1}{2} \, c^2(u) \, u_x^2 + \frac{1}{2} \, b^2(u) \, u_y^2 + a(u) \, u_x \, u_y.
$$ 
Here, $H$ is the ``Hamiltonian'', and $\frac{\delta H}{\delta u}$ denotes the variational
derivative of function $H(u, u_x, u_y)$ with respect to $u$.

A simple calculation, in light of \eqref{eq:vari}, reveals that
\begin{equation}
  \label{eq:energycont}
  \frac{d}{dt} \int_{\R} \left( \frac{1}{2} u_t^2 + H(u, u_x, u_y) \right) \,dx =0.
\end{equation}
To be more precise, this is a direct consequence of the simple identity:
\begin{equation}
  \label{eq:varderivative}
  \frac{\delta H}{\delta u} = \frac{\partial H}{\partial u} - \frac{d}{dx}\left(\frac{\partial H}{\partial u_x}\right) - \frac{d}{dy}\left(\frac{\partial H}{\partial u_y}\right).
\end{equation}
We also note that for equation \eqref{eq:main_2d},
\begin{align*}
  \frac{\delta H}{\delta u} & = c(u) c'(u) u_x^2 - \left( c^2(u)
    u_x\right)_x + b(u) b'(u) u_y^2 - \left( b^2(u)
    u_y\right)_y + a'(u) u_x u_y -(a(u) u_y)_x - (a(u) u_x)_y  \\
   & =- c^2(u) u_{xx} - c(u) c'(u) u_x^2  - b^2(u) u_{yy} - b(u) b'(u) u_y^2 - a'(u) u_x u_y - 2 a(u) u_{xy} \\
   &= - c(u) \left( c(u) u_x \right)_x - b(u) \left( b(u) u_y \right)_y - a'(u) u_x u_y - 2 a(u) u_{xy}.
\end{align*}
Based on above observations, we propose the following scheme for
\eqref{eq:main_2d}
\begin{equation} 
\begin{aligned}
(u_{ij})_{tt} & + c(u_{ij}) c'(u_{ij}) (D^x u_{ij})^2 - D^x \left(c^2(u_{ij}) D^x u_{ij} \right) + b(u_{ij}) b'(u_{ij}) (D^y u_{ij})^2 - D^y \left(b^2(u_{ij}) D^y u_{ij} \right) \\
&\qquad + a'(u_{ij}) D^x(u_{ij}) D^y(u_{ij}) - D^x \left( a(u_{ij}) D^y u_{ij} \right) - D^y \left( a(u_{ij}) D^x u_{ij} \right) = 0,
\end{aligned}
\label{eq:scheme9}
\end{equation}
where the central differences $D^x$ and $D^y$ are defined by 
\begin{equation*}
D^x z_{ij} = \frac{z_{i+1,j}-z_{i-1,j}}{2\Dx}, \, \, \text{and} \, \, D^y z_{ij} = \frac{z_{i,j+1}-z_{i,j-1}}{2\Dy}.
\end{equation*}
This scheme is  energy preserving as shown in the following
theorem:
\begin{theorem}
Let $u_{ij}(t)$ be an approximate solution generated by the scheme
\eqref{eq:scheme9} using periodic boundary conditions. Then we have
\begin{align*}
\frac{d}{dt} \left( \frac{\Dx \Dy}{2} \sum_{i,j} (u_{ij})_t^2 + c^2(u_{ij}) \left( D^x u_{ij} \right)^2  + b^2(u_{ij}) \left( D^y u_{ij} \right)^2 + 2 a(u_{ij}) D^x (u_{ij}) D^y (u_{ij}) \right) = 0.
\end{align*}
\end{theorem}

\begin{proof}
We start by calculating
\begin{align*}
\frac{d}{dt} & \left( \frac{\Dx \Dy}{2} \sum_{i,j} (u_{ij})_t^2 + c^2(u_{ij}) \left( D^x u_{ij} \right)^2  + b^2(u_{ij}) \left( D^y u_{ij} \right)^2 + 2 a(u_{ij}) D^x (u_{ij}) D^y (u_{ij}) \right) \\
&= \Dx \Dy \sum_{i,j} \left( (u_{ij})_t (u_{ij})_{tt}  + c(u_{ij}) c'(u_{ij}) \left( D^x u_{ij} \right)^2 (u_{ij})_t + c^2(u_{ij}) D^x u_{ij} D^x (u_{ij})_t \right) \\
& \qquad + \Dx \Dy \sum_{i,j} \left( b(u_{ij}) b'(u_{ij}) \left( D^y u_{ij} \right)^2 (u_{ij})_t + b^2(u_{ij}) D^y u_{ij} D^y (u_{ij})_t \right) \\
& \, +  \Dx \Dy \sum_{i,j} \left( a'(u_{ij}) D^x(u_{ij}) D^y(u_{ij}) (u_{ij})_t  +  a(u_{ij}) D^x (u_{ij})_t  D^y u_{ij} + a(u_{ij}) D^x u_{ij} D^y (u_{ij})_t   \right) \\
& = \Dx \Dy \sum_{i,j} \left( (u_{ij})_t (u_{ij})_{tt} + c(u_{ij}) c'(u_{ij}) \left( D^x u_{ij} \right)^2 (u_{ij})_t - D^x \left(c^2(u_{ij}) D^x u_{ij}\right) (u_{ij})_t \right) \\
& \qquad + \Dx \Dy \sum_{i,j} \left( b(u_{ij}) b'(u_{ij}) \left( D^y u_{ij} \right)^2 (u_{ij})_t -  D^y \left(b^2(u_{ij}) D^y u_{ij}\right) (u_{ij})_t \right) \\    
& \, + \Dx \Dy \sum_{i,j} \left( a'(u_{ij}) D^x(u_{ij}) D^y(u_{ij}) (u_{ij})_t  -  D^x \left( a(u_{ij})  D^y u_{ij} \right) (u_{ij})_t - D^y \left( a(u_{ij}) D^x u_{ij} \right) (u_{ij})_t   \right) \\
& = 0. \, \, \text{(follows from \eqref{eq:scheme9})}
  \end{align*}

\end{proof}

\section{Numerical Experiments} 
\label{sec:numerical}
For the numerical experiments, the computational domain is subdivided
into $N \times N$ rectangular cells. All cells are of size $\Delta x
\times \Delta y$. A uniform time step 
\begin{equation}
  \Delta t = 0.1 \frac{\min\{\Delta x,\Delta y\}}{\max\{\alpha,\beta\}}
\end{equation}
is used throughout the computation. Moreover, in all experiments
the parameters for the shock capturing operator are $C = 0.1$ and $\theta = 1$.
To keep focus on the spatial discretization, we will use a fifth-order
Runge--Kutta scheme \cite{Lut} ensuring a satisfactory temporal accuracy.
Periodic boundary conditions are used in all experiments.

\subsection{Gaussian disturbance to homogeneous director state}
\begin{figure}[htbp]
  \centering
  \begin{subfigure}[b]{0.45\textwidth}
    \includegraphics[width=\textwidth]{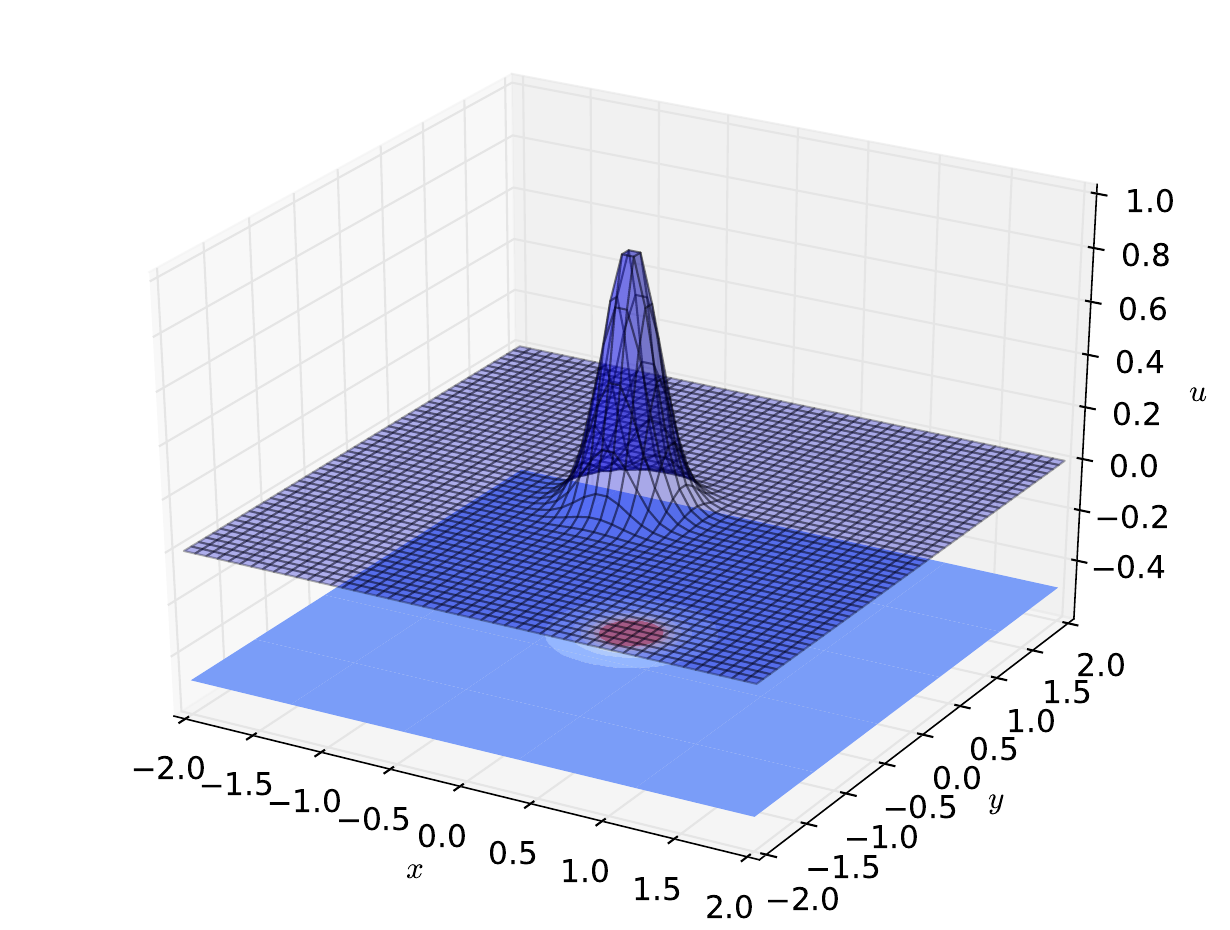}
    \caption{$t = 0$}
  \end{subfigure}
  \begin{subfigure}[b]{0.45\textwidth}
    \includegraphics[width=\textwidth]{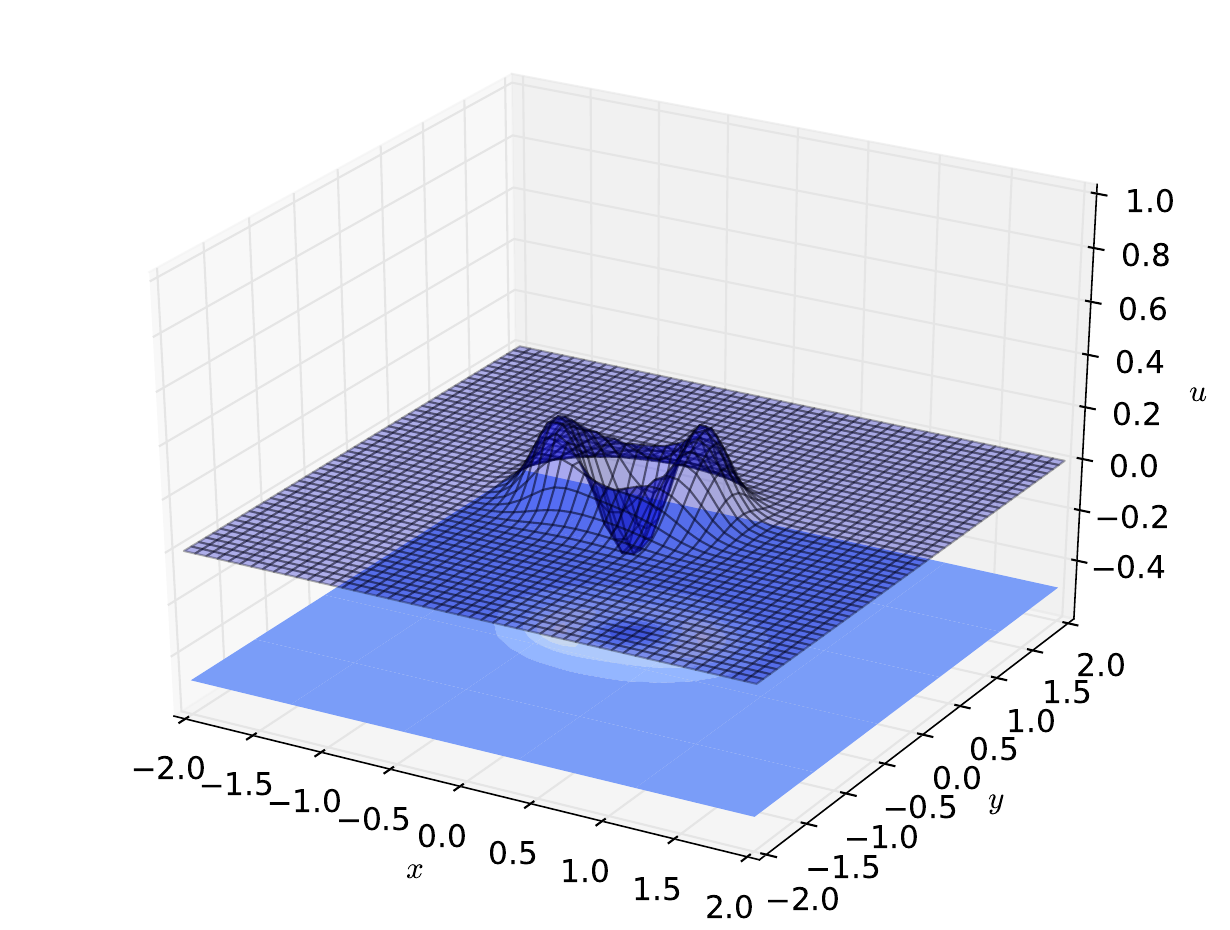}
    \caption{$t = 1/3$}
  \end{subfigure}
  \begin{subfigure}[b]{0.45\textwidth}
    \includegraphics[width=\textwidth]{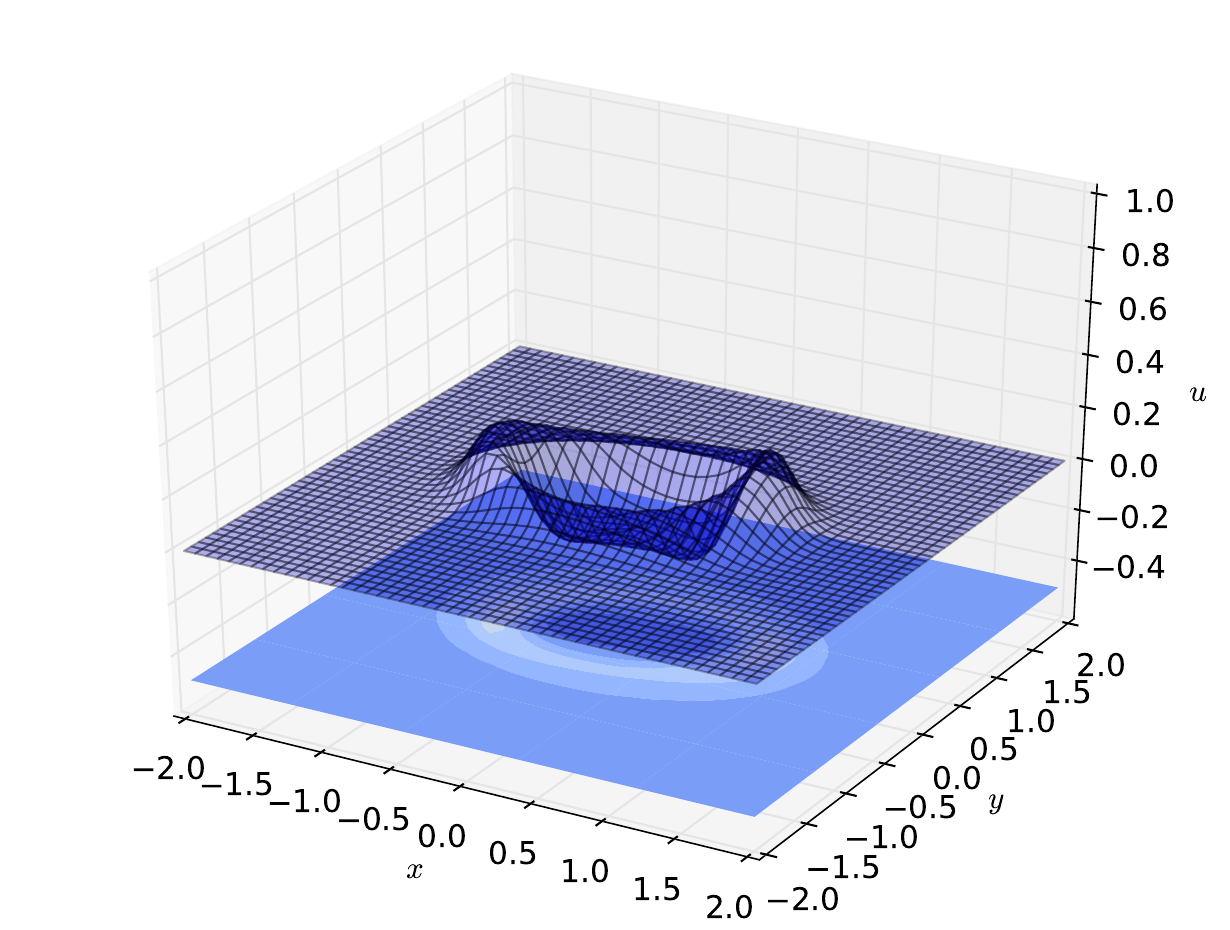}
    \caption{$t = 2/3$}
  \end{subfigure}
  \begin{subfigure}[b]{0.45\textwidth}
    \includegraphics[width=\textwidth]{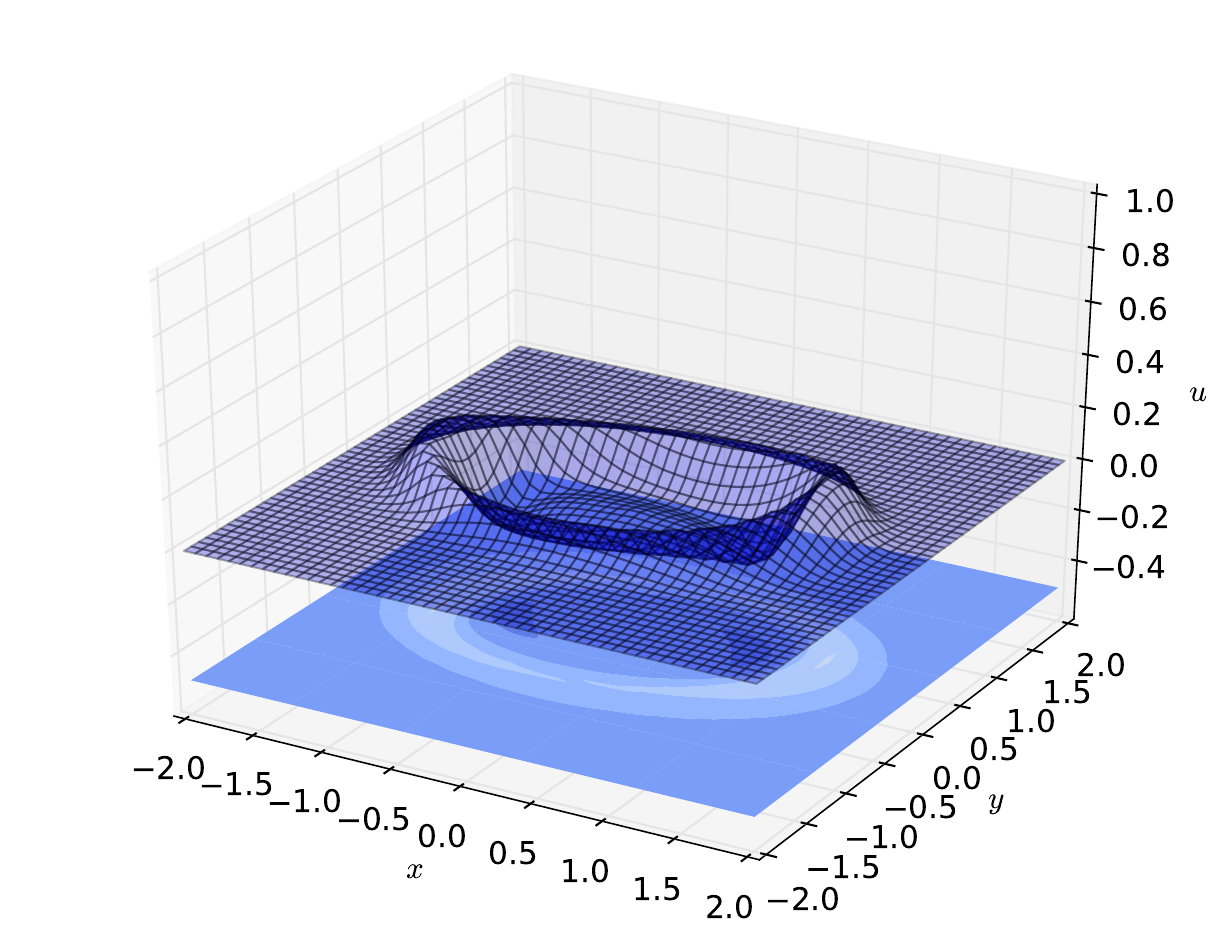}
    \caption{$t = 1$}
  \end{subfigure}
  \caption{Numerical solution of the initial value problem \eqref{eq:main_2d}
    with the initial data \eqref{eq:gaussianData} using the
    dissipative piecewise quadratic scheme with $N = 32$. 
    The parameters are $\alpha = 1.5$ and $\beta = 0.5$.}
  \label{fig:gaussian_d2}
\end{figure}

In this section we consider the initial value problem \eqref{eq:main_2d} with
the initial data
\begin{subequations}
\begin{align}
  u_0(x,y) &= \exp\left(-16 \left(x^2 + y^2\right) \right)\\
  u_1(x,y) &= 0
\end{align}
\label{eq:gaussianData}
\end{subequations}
on $(x,y) \in \mathbb{R}^2$. The physical parameters are $\alpha = 1.5$ and
$\beta = 0.5$. A numerical solution was computed using $N = 32$ with
the dissipative piecewise quadratic ($s = 2$) scheme. Figure \ref{fig:gaussian_d2} 
shows the time evolution of the numerical solution, demonstrating the
non-isotropic nature of this model. 

\begin{figure}[htbp]
  \centering
    \includegraphics[width=0.45\textwidth]{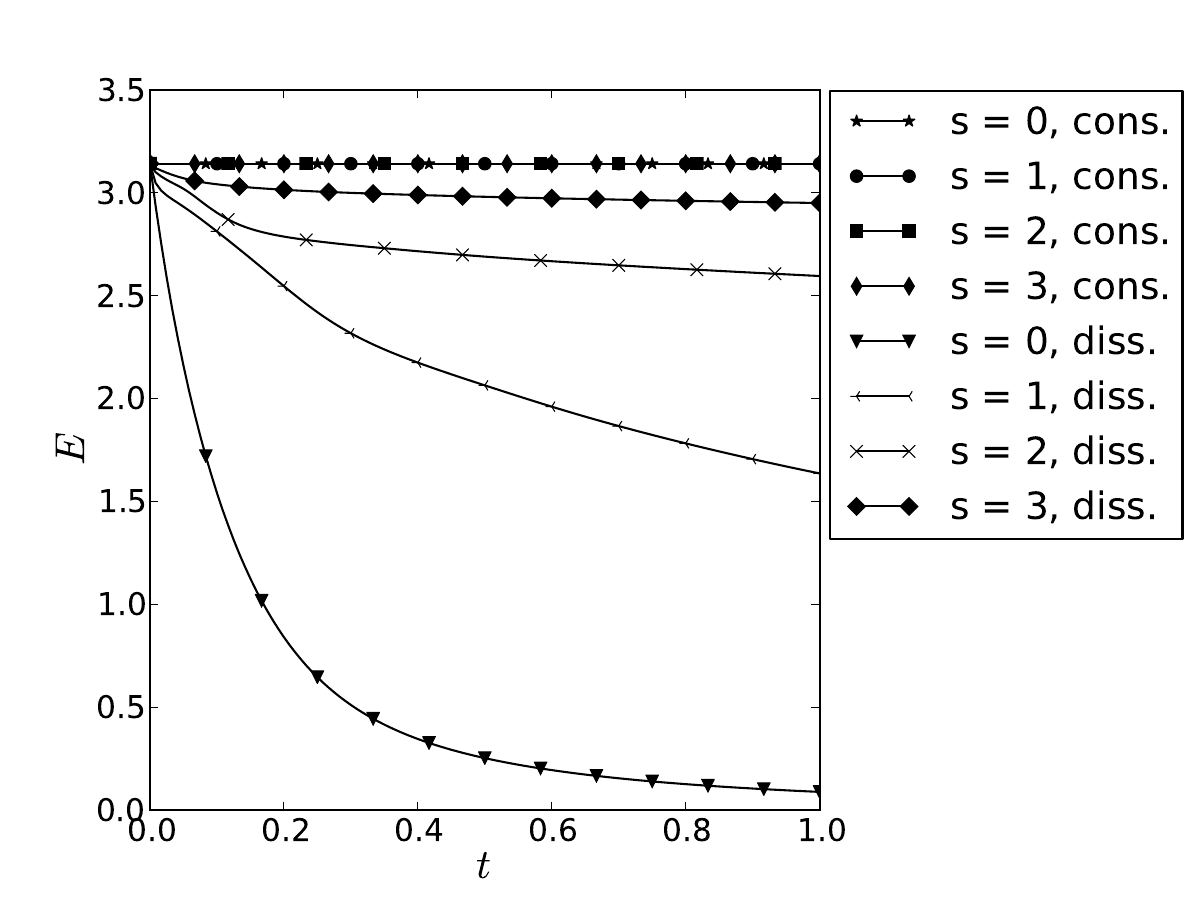}
    \caption{Evolution of the discrete energy \eqref{eq:discrete_energy}
    for the numerical solutions of the 
    initial value problem \eqref{eq:main_2d} with the initial data
    \eqref{eq:gaussianData} using both conservative and dissipative schemes.
    The parameters were
    $\alpha = 1.5$ and $\beta = 0.5$ and a $N = 32$ grid size was used.}
  \label{fig:gaussianEnergy}
\end{figure}

A key property of the schemes derived in this paper is that they are designed,
at the semi-discrete level, to either conserve or dissipate the energy. 
Figure \ref{fig:gaussianEnergy} shows the time evolution of the discrete
energy
\begin{equation} \label{eq:discrete_energy}
  E = \sum_{i,j=1}^N \int_{\Omega_{ij}} \frac{p^2 + \alpha v^2
  + \beta w^2}{2} \mathrm{d} x
  = \frac{\Delta x \, \Delta y}{8} \sum_{i,j=1}^N \sum_{k,l=0}^s \rho_k 
  \rho_l \left( \left(p_{ij}^{(kl)}\right)^2 + \alpha \left(v_{ij}^{(kl)}
     \right)^2
     + \beta \left(w_{ij}^{(kl)}\right)^2 \right),
\end{equation}
for the Gaussian initial value problem using both conservative and
dissipative schemes for ${s \in \{0,\cdots,3 \}}$. The results clearly indicate
that the energy preserving (and dissipating) properties carry over to the fully
discrete case when using a higher-order time integrator.

\subsection{Loss of regularity} 
\begin{figure}[htbp]
  \captionsetup[subfigure]{labelformat=empty}
  \centering
  \begin{subfigure}{0.90\textwidth}
  \begin{subfigure}[b]{0.5\textwidth}
    \includegraphics[width=\textwidth]{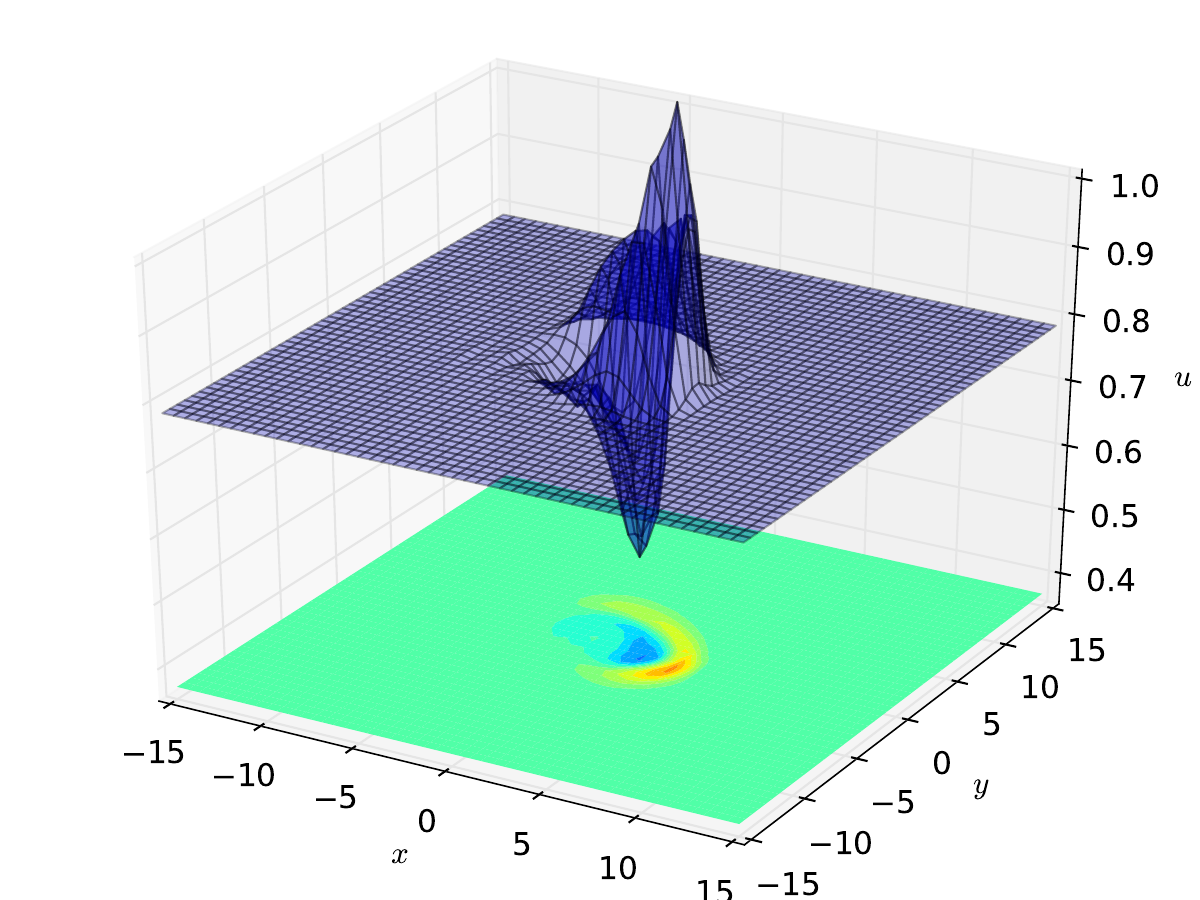}
    \caption{Conservative scheme}
  \end{subfigure}
  \begin{subfigure}[b]{0.5\textwidth}
    \includegraphics[width=\textwidth]{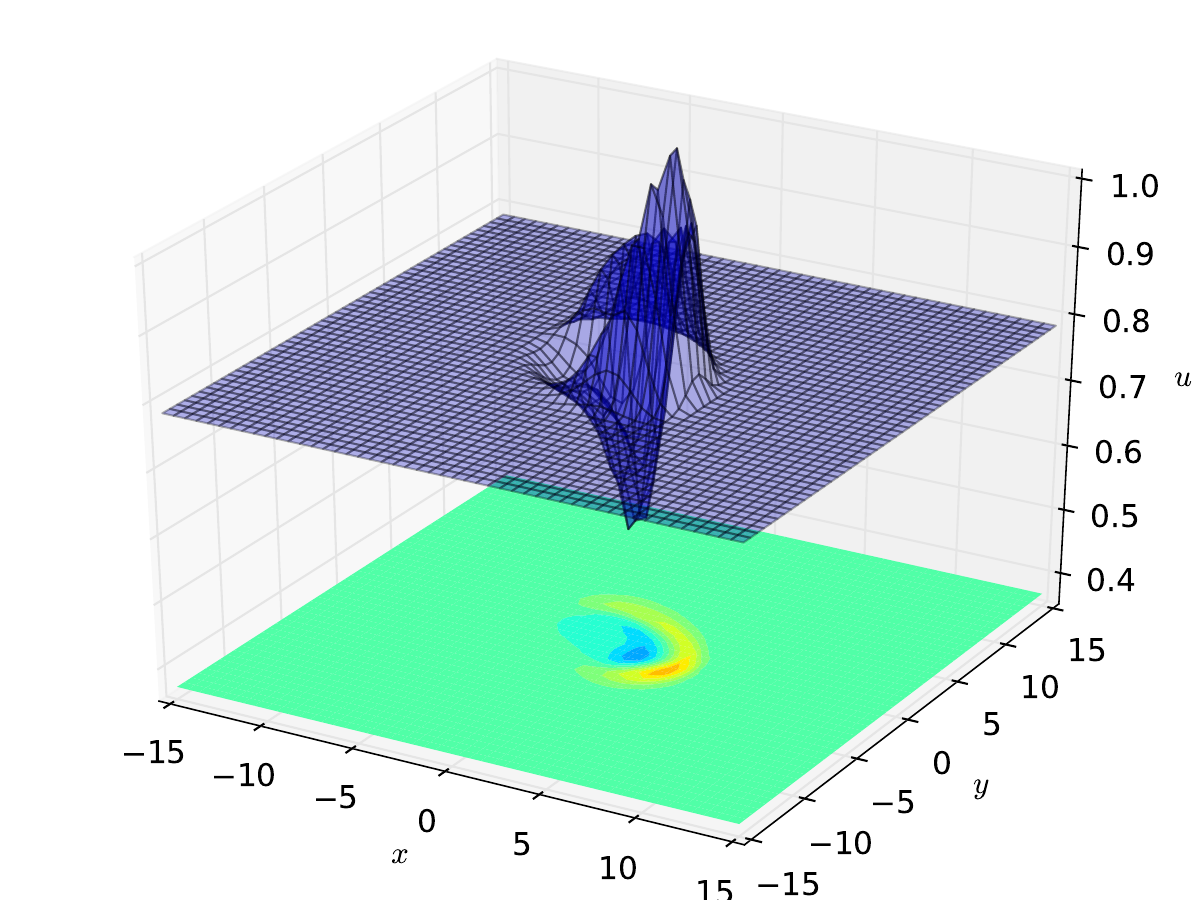}
    \caption{Dissipative scheme}
  \end{subfigure}
  \caption{$t = 3.33$}
  \end{subfigure}
  \begin{subfigure}[b]{0.90\textwidth}
  \begin{subfigure}[b]{0.5\textwidth}
    \includegraphics[width=\textwidth]{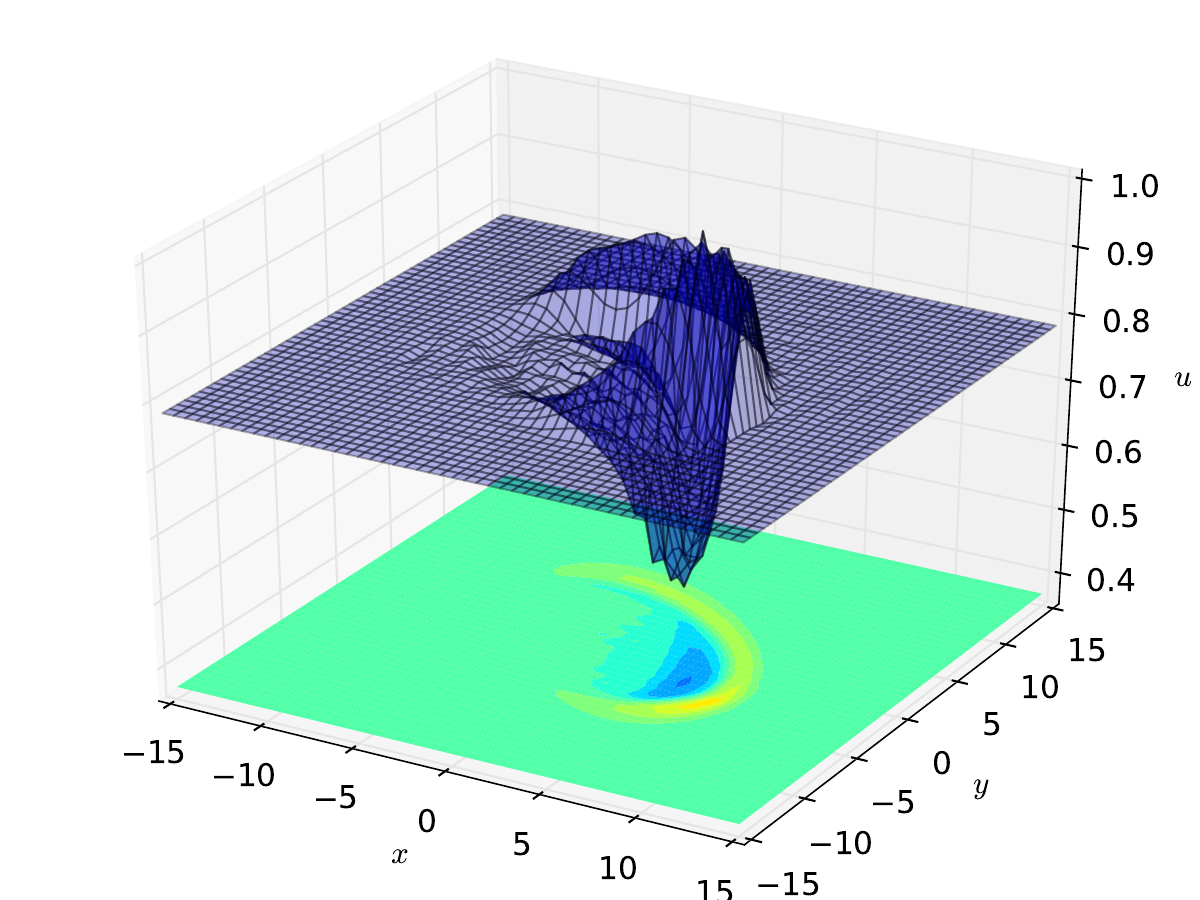}
    \caption{Conservative scheme}
  \end{subfigure}
  \begin{subfigure}[b]{0.5\textwidth}
    \includegraphics[width=\textwidth]{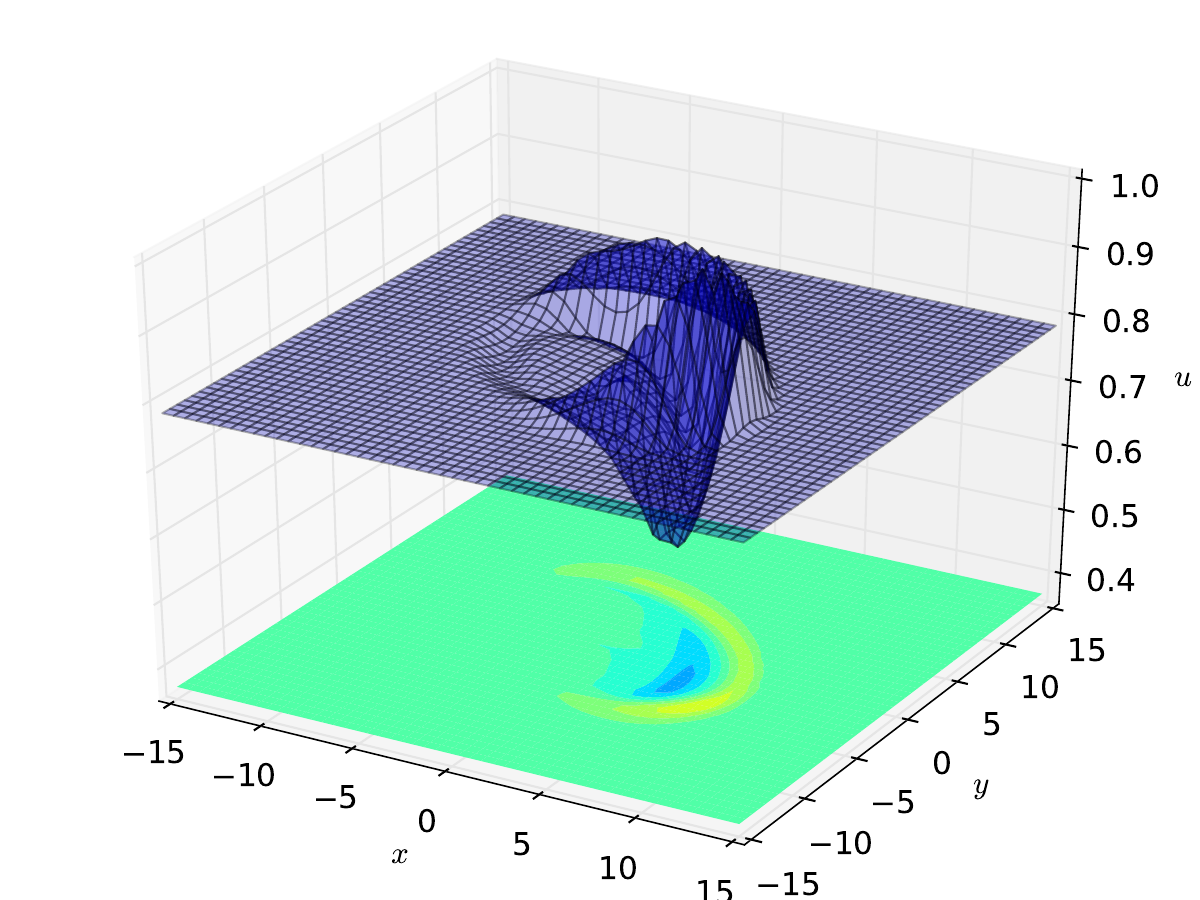}
    \caption{Dissipative scheme}
  \end{subfigure}
  \caption{$t = 6.66$}
  \end{subfigure}
  \begin{subfigure}[b]{0.90\textwidth}
  \begin{subfigure}[b]{0.5\textwidth}
    \includegraphics[width=\textwidth]{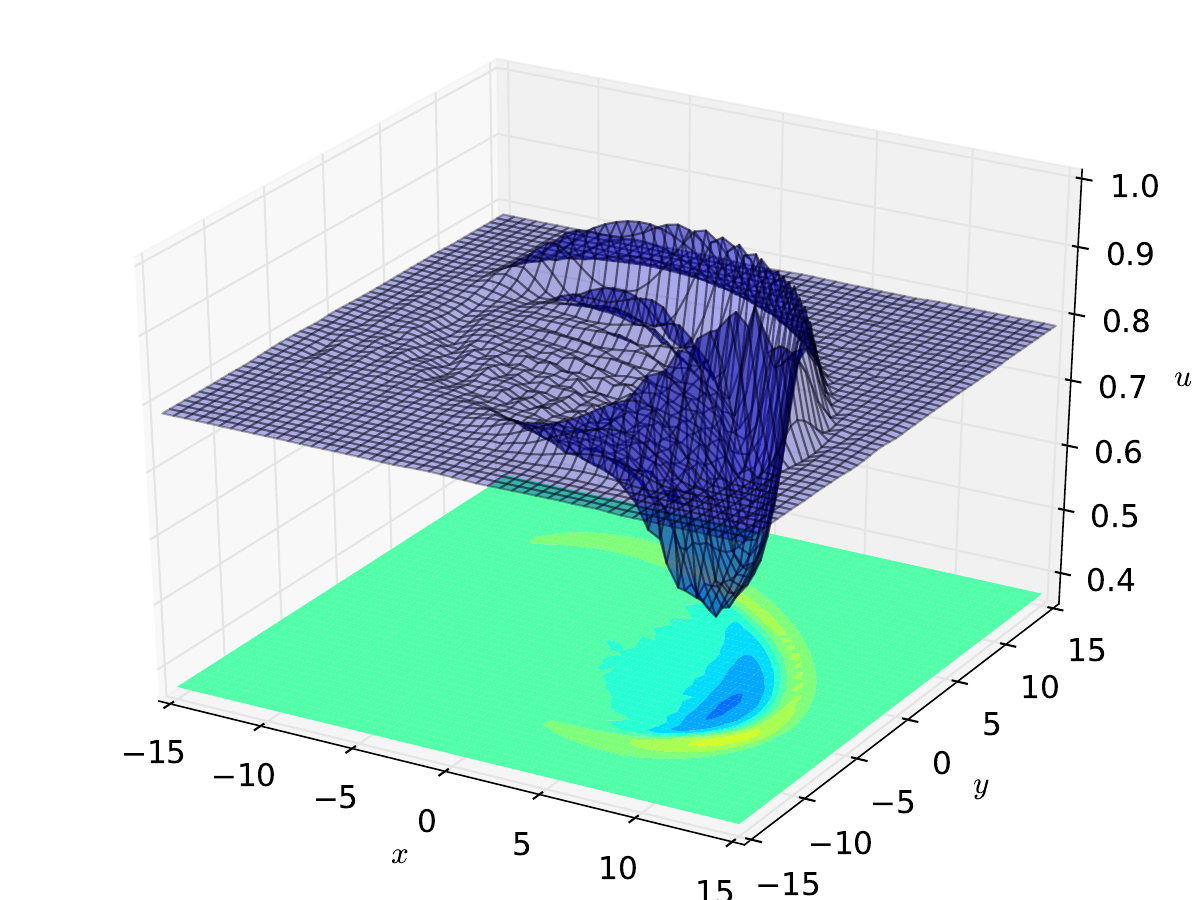}
    \caption{Conservative scheme}
  \end{subfigure}
  \begin{subfigure}[b]{0.5\textwidth}
    \includegraphics[width=\textwidth]{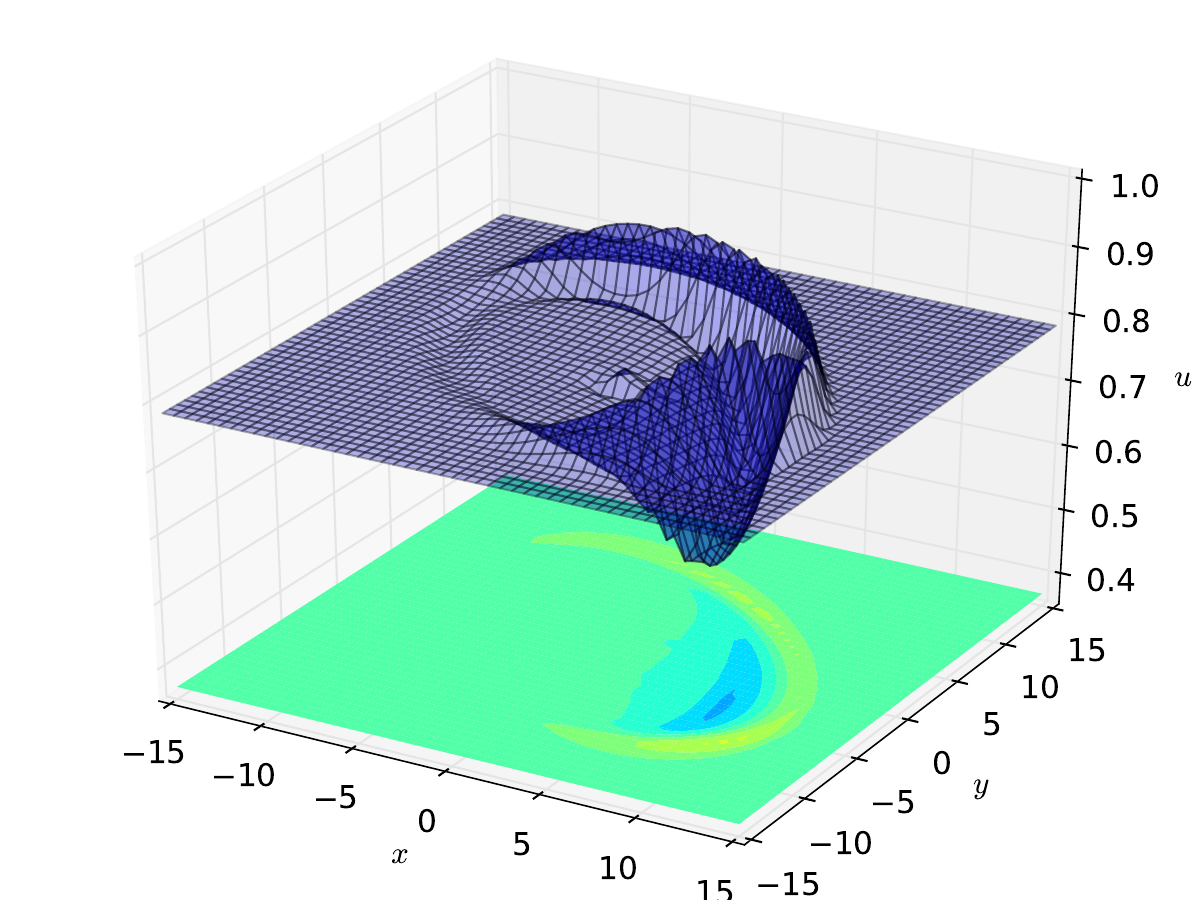}
    \caption{Dissipative scheme}
  \end{subfigure}
  \caption{$t = 10$}
  \end{subfigure}
  \caption{The numerical solution at for left) the piecewise quadratic
    conservative scheme and right) the piecewise quadratic dissipative scheme 
    of the initial value problem
    \eqref{eq:main_2d} with initial data \eqref{eq:gaussianTravellingData} 
    with $N = 64$ cells. 
    The physical parameters were $\alpha = 1.5$ and $\beta = 0.5$.}
  \label{fig:gaussian_travelling_case}
\end{figure}

A crucial property for the 1D variational wave equation is that solutions
loose regularity in finite time even for smooth initial data. For the 2D case
this is still an open problem. We investigate this numerically by considering
the initial value problem \eqref{eq:main_2d} with data 
\begin{subequations}
\begin{align}
  u_0(x,y) &= \exp\left(-\left(x^2 + y^2\right) \right)\\
  u_1(x,y) &= -c(u_0(x,y)) u_{0,x}(x,y)
\end{align}
\label{eq:gaussianTravellingData}
\end{subequations}
for $(x,y) \in \mathbb{R}^2$. A numerical experiment was performed using
$N = 64$ computational cells with the conservative and dissipative piecewise
quadratic schemes. The results, shown in Figure
\ref{fig:gaussian_travelling_case}, indicates a clear steepening of the gradient
as the solution evolves. 

\begin{figure}[htbp]
  \centering
  \begin{subfigure}[b]{0.45\textwidth}
    \includegraphics[width=\textwidth]{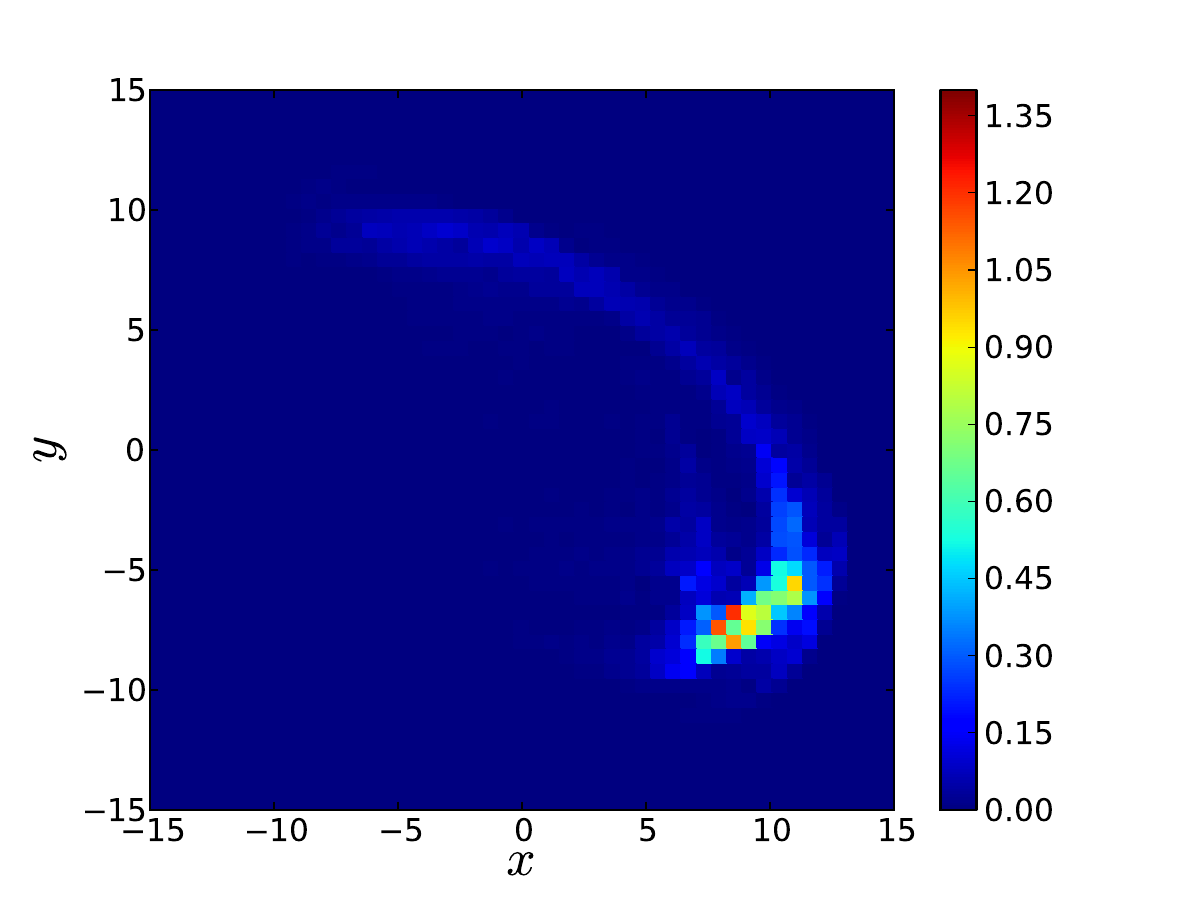}
    \caption{Conservative scheme}
  \end{subfigure}
  \begin{subfigure}[b]{0.45\textwidth}
    \includegraphics[width=\textwidth]{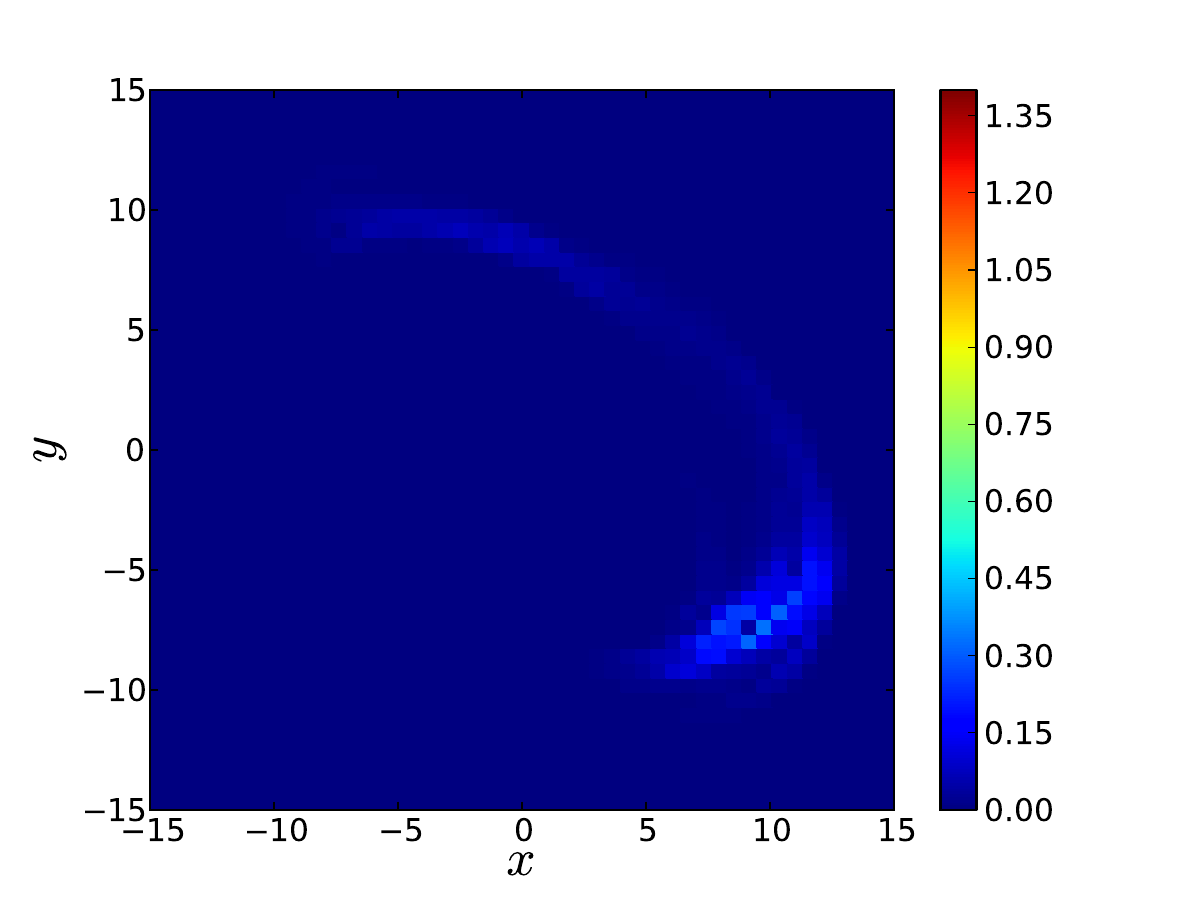}
    \caption{Dissipative scheme}
  \end{subfigure}
  \caption{The root-mean-square of the residual \eqref{eq:residual_rms} at $t =
    10$ for the initial value problem \eqref{eq:main_2d} with initial data 
    \eqref{eq:gaussianTravellingData}. At the left: the piecewise quadratic
    conservative scheme and at the right: the piecewise quadratic dissipative
    scheme, both with $N = 64$ cells. 
    The physical parameters were $\alpha = 1.5$ and $\beta = 0.5$.}
  \label{fig:travelling_residual}
\end{figure}

Smooth solutions of \eqref{eq:main_2d} satisfies the conservation law 
\eqref{eq:residual}. The root-mean-square of the residual
\eqref{eq:residual_rms} can therefore be an indicator function for loss of
regularity in the solution. Figure \ref{fig:travelling_residual} shows the
residual at $t = 10$ for both the conservative and dissipative schemes. The
results indicate that the solution looses smoothness near the front of the
wave propagating in the positive $x$ direction. Moreover, as expected,
the dissipative scheme with
the shock capturing operator is able to maintain a higher degree of numerical
smoothness (as measured by the residual) than the conservative scheme.

\subsection{Bifurcation of solutions}
Another critical feature of the 1D nonlinear variational wave equation
\eqref{eq:main} is the existence of different classes of weak solutions.
However, the existence and well-posedness for the initial value problem in the 
2D generalization remains an open problem. 

\begin{figure}[htbp]
  \centering
  \begin{subfigure}[b]{0.45\textwidth}
    \includegraphics[width=\textwidth]{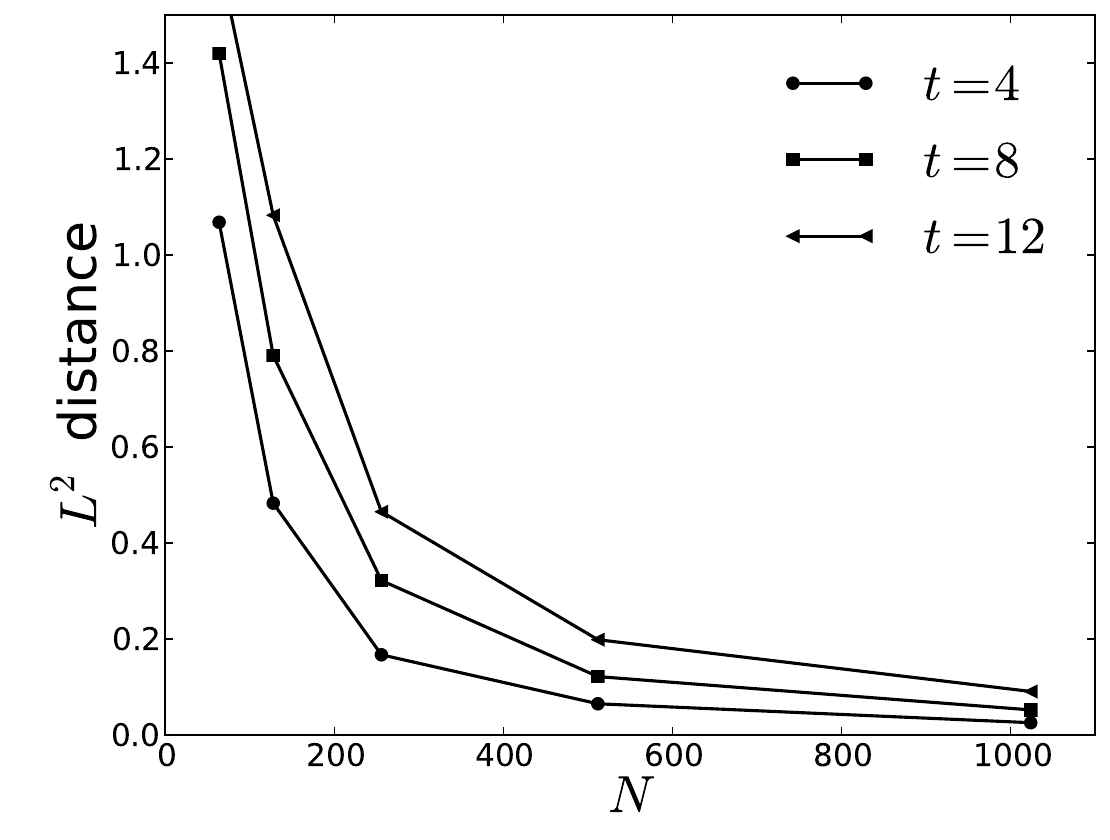}
    \caption{Conservative - Hamiltonian}
  \end{subfigure}
  \begin{subfigure}[b]{0.45\textwidth}
    \includegraphics[width=\textwidth]{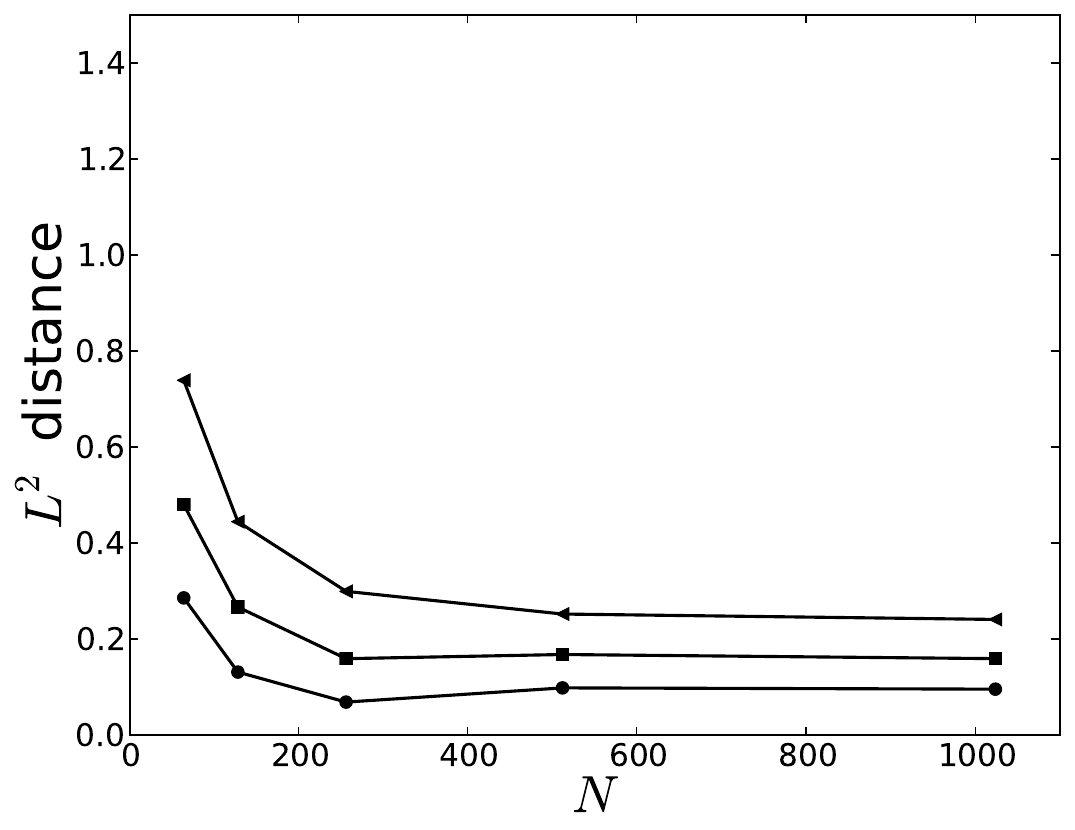}
    \caption{Conservative - dissipative}
  \end{subfigure}
  \caption{The $L^2$ distance between left: the conservative DG scheme and the
    Hamiltonian scheme and right: the conservative DG scheme and the
    dissipative DG scheme, for the initial value problem \eqref{eq:main_2d} 
    with initial data \eqref{eq:gaussianTravellingData}. 
    The physical parameters were $\alpha = 1.5$ and $\beta = 0.5$.}
  \label{fig:travelling_l2dist}
\end{figure}

In order to investigate this issue numerically, we consider the initial data
\ref{fig:gaussian_travelling_case} and study the convergence of the three
schemes; the conservative DG scheme, the dissipative DG scheme and the Hamiltonian 
scheme; after the loss of regularity. Figure \ref{fig:travelling_l2dist} shows
the $L^2$ distance between the numerical solutions for different times and
under grid refinement. The results indicate that the conservative DG scheme
and the Hamiltonian scheme indeed converge to the same solution as the grid is
refined. However, the distance between the dissipative and conservative DG
schemes seems to converge to a non-zero value that increases as a function of
time. This may indicate that the question of well-posedness for the 2D
variational wave equation is as delicate as in the 1D case.

\subsection{Order of Convergence and Efficiency}
\begin{figure}[htbp]
  \centering
  \begin{subfigure}[b]{0.45\textwidth}
    \includegraphics[width=\textwidth]{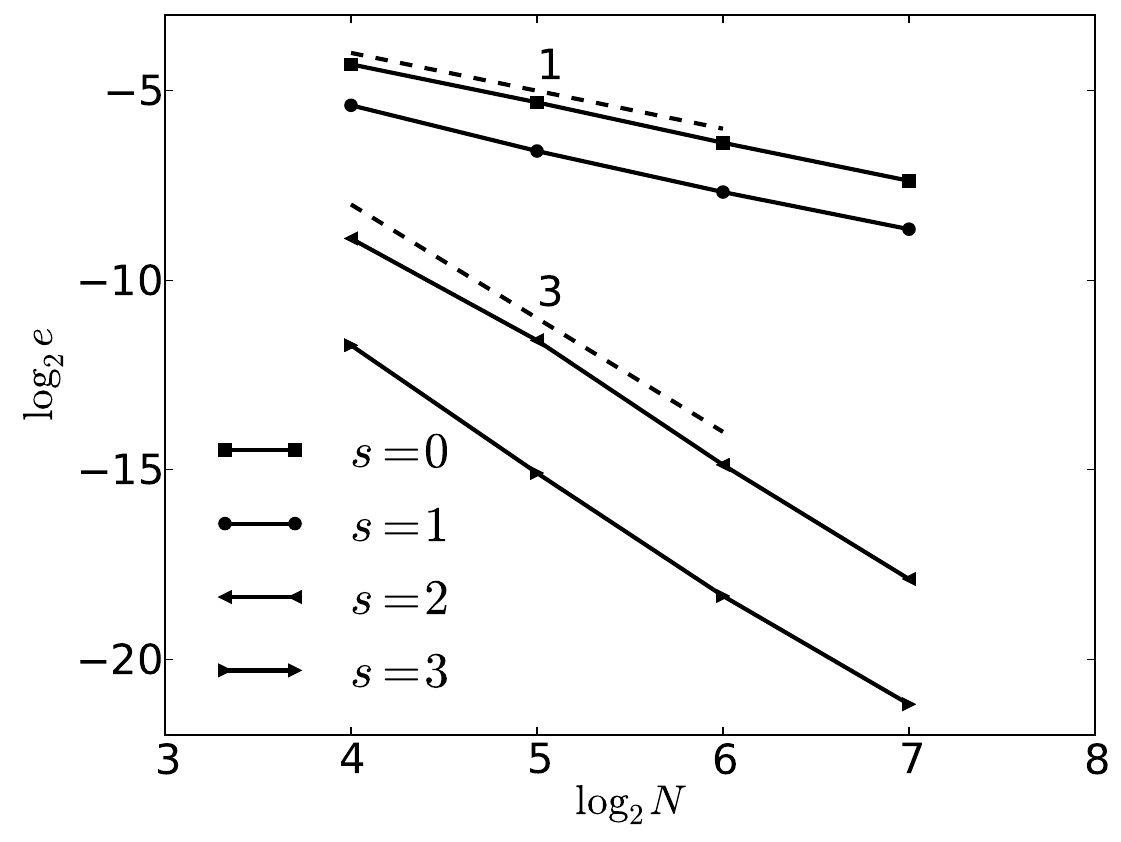}
    \caption{Conservative scheme}
  \end{subfigure}
  \begin{subfigure}[b]{0.45\textwidth}
    \includegraphics[width=\textwidth]{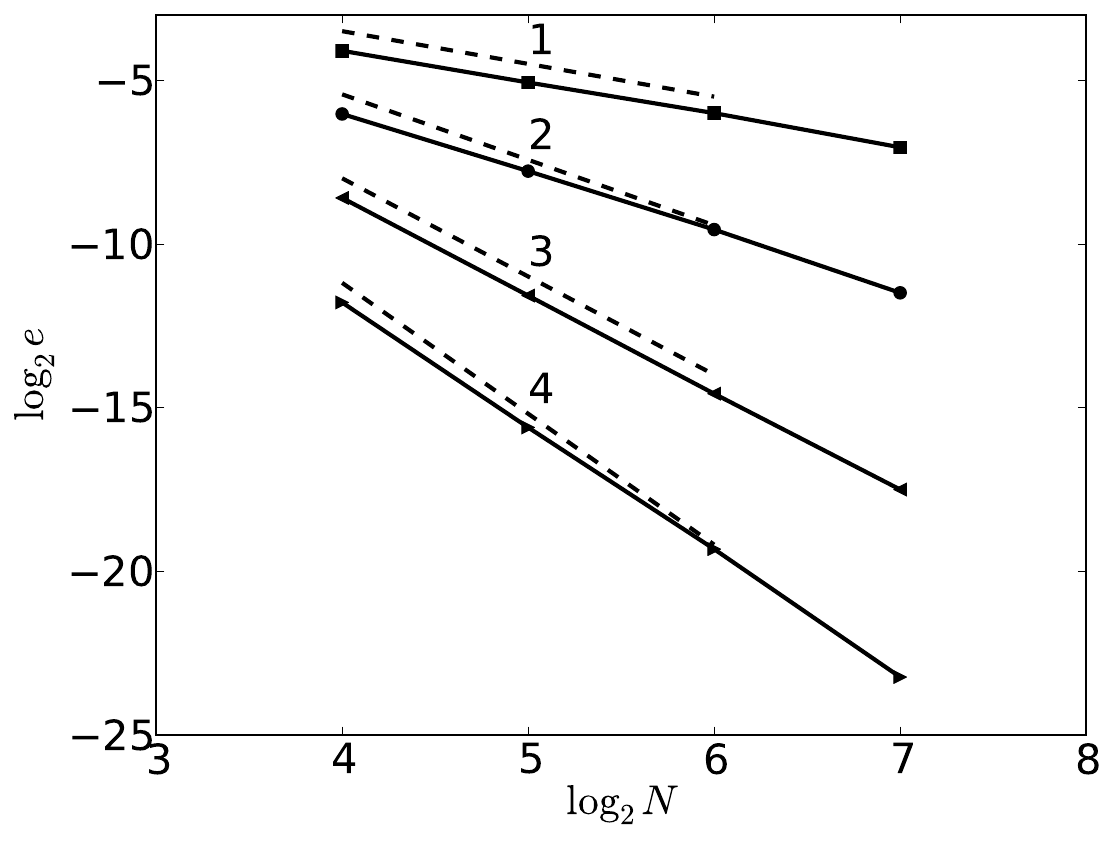}
    \caption{Dissipative scheme}
  \end{subfigure}
  \begin{subfigure}[b]{0.45\textwidth}
    \includegraphics[width=\textwidth]{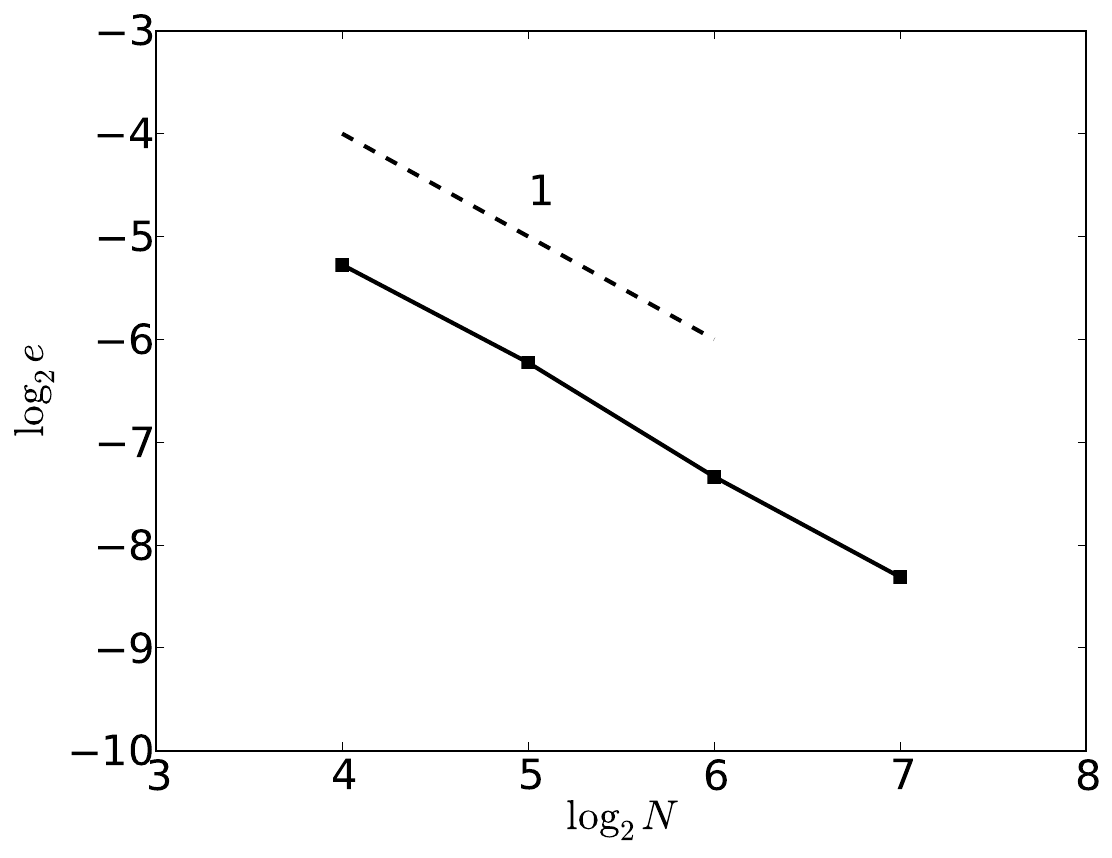}
    \caption{Hamiltonian scheme}
  \end{subfigure}
  \caption{The error \eqref{eq:error} for the numerical solution of the 
  Gaussian initial value
  problem \eqref{eq:gaussianData} 
  as a function of $N$, using $\alpha = 1.5$ and 
  $\beta = 0.5$. The dashed lines indicate the different orders of 
  convergence.}
  \label{fig:orderOfConvergence}
\end{figure}

In the following, we demonstrate the order of convergence and efficiency 
of both the conservative and dissipative schemes for smooth solutions. 
As before, we consider the initial value problem \eqref{eq:main_2d} with
the initial data \eqref{eq:gaussianData} 
with physical parameters $\alpha = 1.5$ and $\beta = 0.5$. A reference
solution $u_\text{ref}$ was calculated at $t = 0.1$ using the conservative piecewise 
cubic scheme ($s = 3$) with $N = 1024$. 
Figure \ref{fig:orderOfConvergence} shows the error
\begin{equation}
  e = \| u_{N} - u_\text{ref} \|_2
  \label{eq:error}
\end{equation}
for different grid cell numbers $N = N_x = N_y$.
The results indicate a suboptimal order of convergence for odd $s$ when using
the conservative numerical flux. For the dissipative scheme the order of
convergence is optimal. This behavior has been observed also in the 1D case
\cite{Aursand2014Preprint}, and for certain DG schemes in the literature
\cite{shu}. The Hamiltonian scheme converges to first order.

\begin{figure}[htbp]
  \centering
  \begin{subfigure}[b]{0.45\textwidth}
    \includegraphics[width=\textwidth]{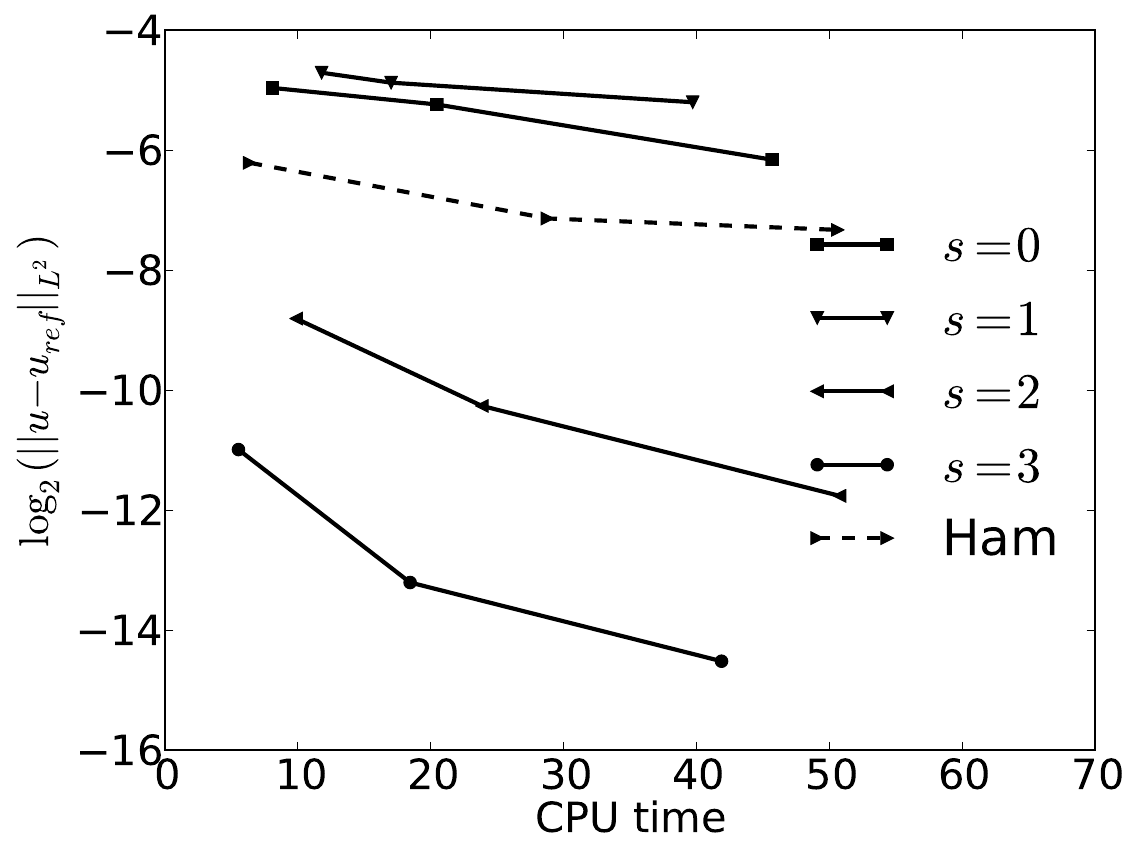}
    \caption{Conservative scheme}
  \end{subfigure}
  \begin{subfigure}[b]{0.45\textwidth}
    \includegraphics[width=\textwidth]{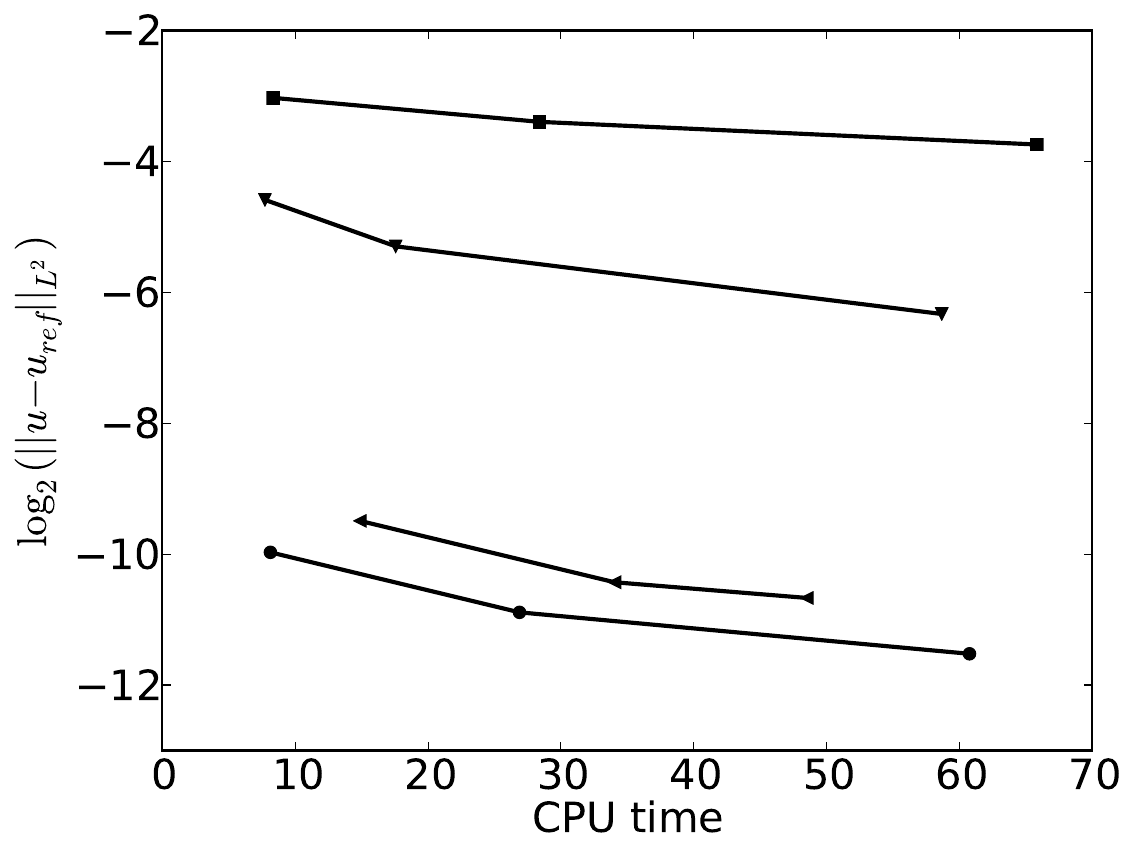}
    \caption{Dissipative scheme}
  \end{subfigure}
  \caption{The error \eqref{eq:error} for the numerical solution of the 
  Gaussian initial value problem \eqref{eq:gaussianData} at $t= 0.5$ as 
  a function of CPU time (wall time), using $\alpha = 1.5$ and 
  $\beta = 0.5$. The reference solution was calculated using the piecewise
  cubic conservative scheme with $N = 1024$ cells.}
  \label{fig:efficiency}
\end{figure}

Figure \ref{fig:efficiency} shows the error \eqref{eq:error} compared to a
a reference solution as a function of computational cost (CPU wall time). The
results indicate that the higher-order schemes mostly make up for their 
increased computational complexity in better accuracy per CPU time. One
exception is the conservative piecewise linear scheme, which for this case
requires more computational work than the piecewise constant scheme in order
to obtain the same accuracy. A possible explanation for this is that enforcing
energy preservation using piecewise linear elements results in an un-physically 
jagged solution in certain regions. This happens despite the fact that the
converged solution does not exhibit this behavior. For the piecewise linear 
dissipative scheme, this effect is suppressed by the added artificial viscosity.

\subsection{Relaxation from a standing wave}
For this experiment we consider the initial value problem 
\begin{subequations}
\begin{align}
  u_0(x,y) &= 2 \cos(2 \pi x) \sin( 2 \pi x) \\
  u_1(x,y) &= \sin(2 \pi (x - y))
\end{align}
\label{eq:sineData}
\end{subequations}
on $(x,y) \in [0,1] \times [0,1]$ with periodic boundary conditions. The
initial value problem can be seen as describing the following: Initially, a
standing wave is induced in the director field using e.g.~an external
electromagnetic field or mechanical vibrations. At $t = 0$, the external
influence is removed, and the evolution of the director is purely governed by
elastic forces. 

\begin{figure}[htbp]
  \captionsetup[subfigure]{labelformat=empty}
  \centering
  \begin{subfigure}{0.90\textwidth}
  \begin{subfigure}[b]{0.5\textwidth}
    \includegraphics[width=\textwidth]{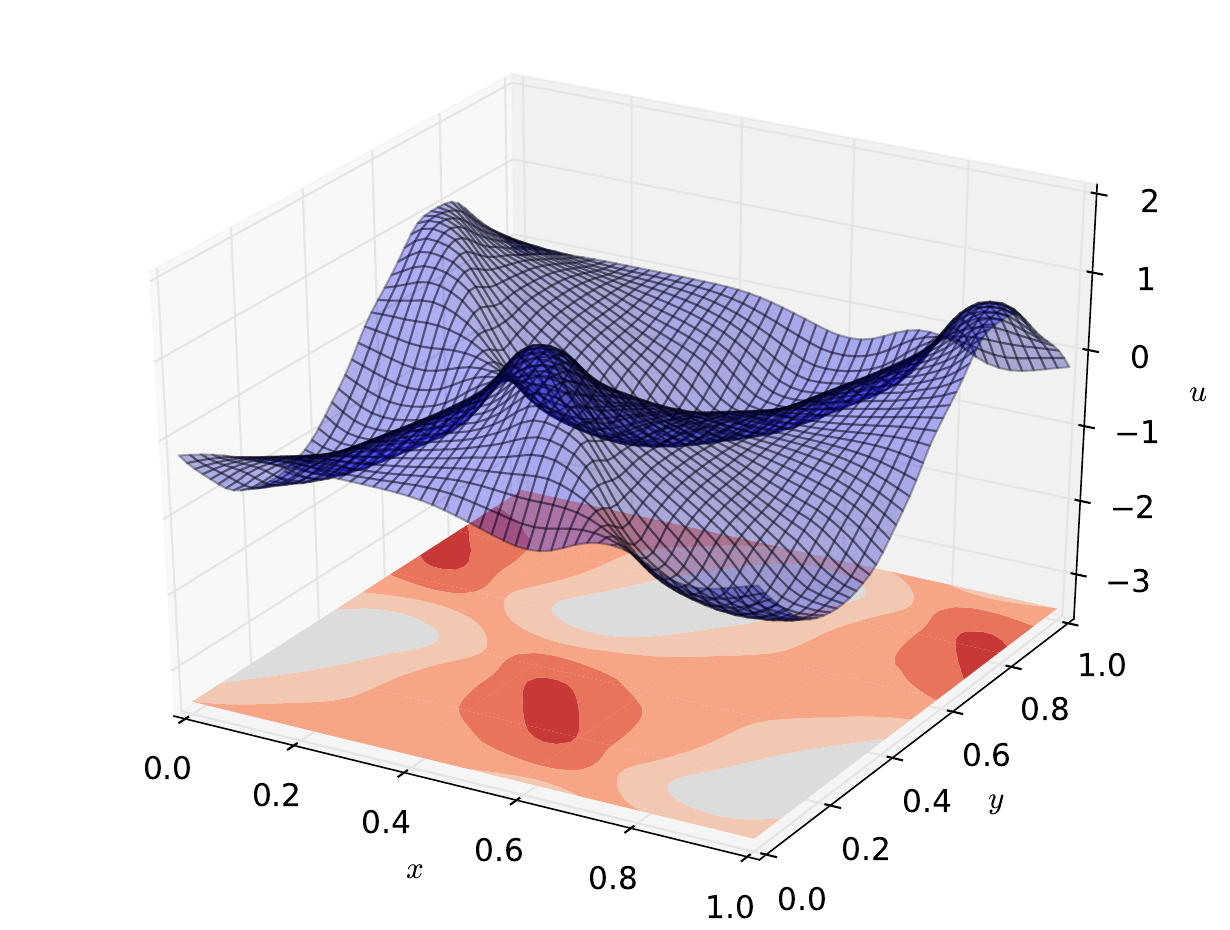}
    \caption{$t = 1$}
  \end{subfigure}
  \begin{subfigure}[b]{0.5\textwidth}
    \includegraphics[width=\textwidth]{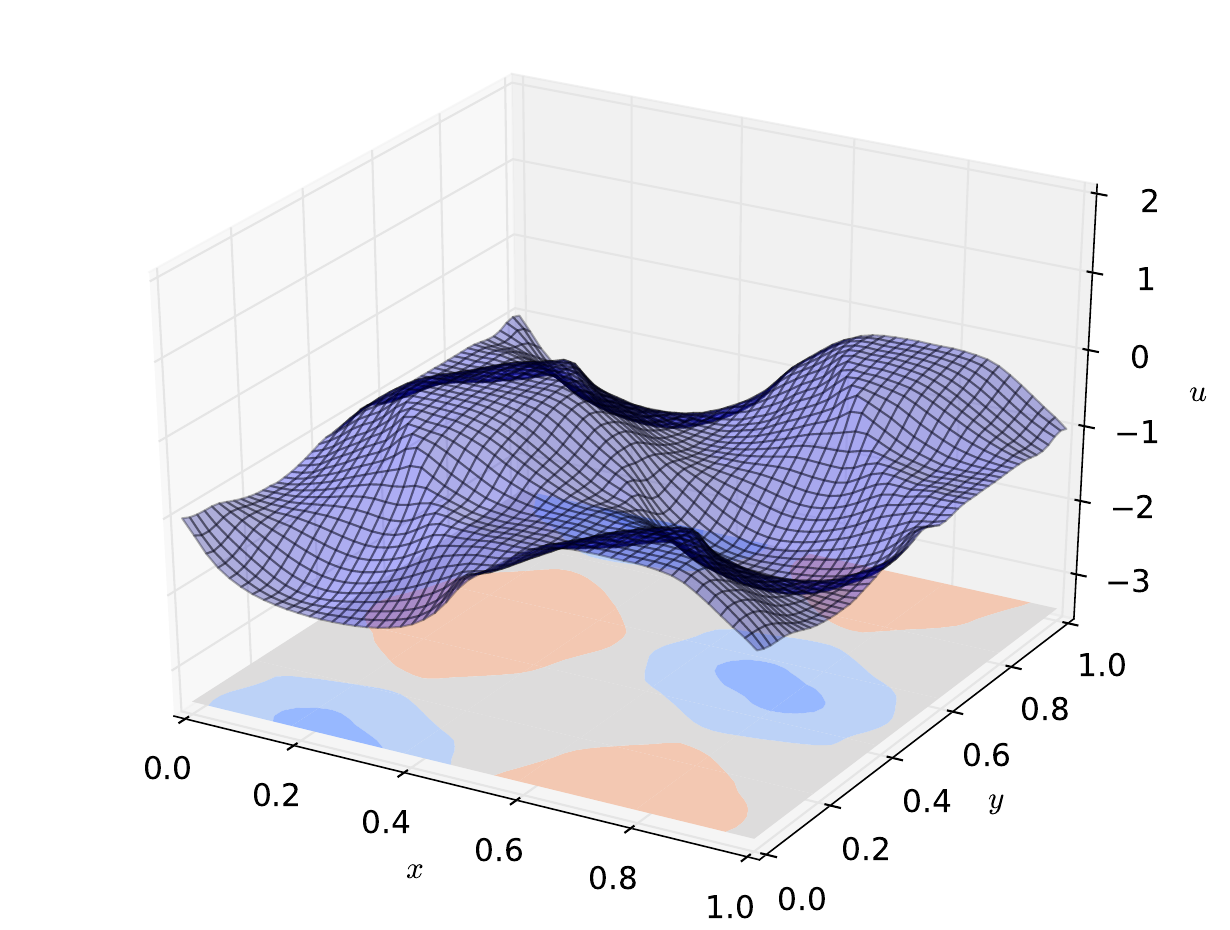}
    \caption{$t = 2$}
  \end{subfigure}
  \caption{Dissipative DG scheme}
  \end{subfigure}
  \begin{subfigure}[b]{0.90\textwidth}
  \begin{subfigure}[b]{0.5\textwidth}
    \includegraphics[width=\textwidth]{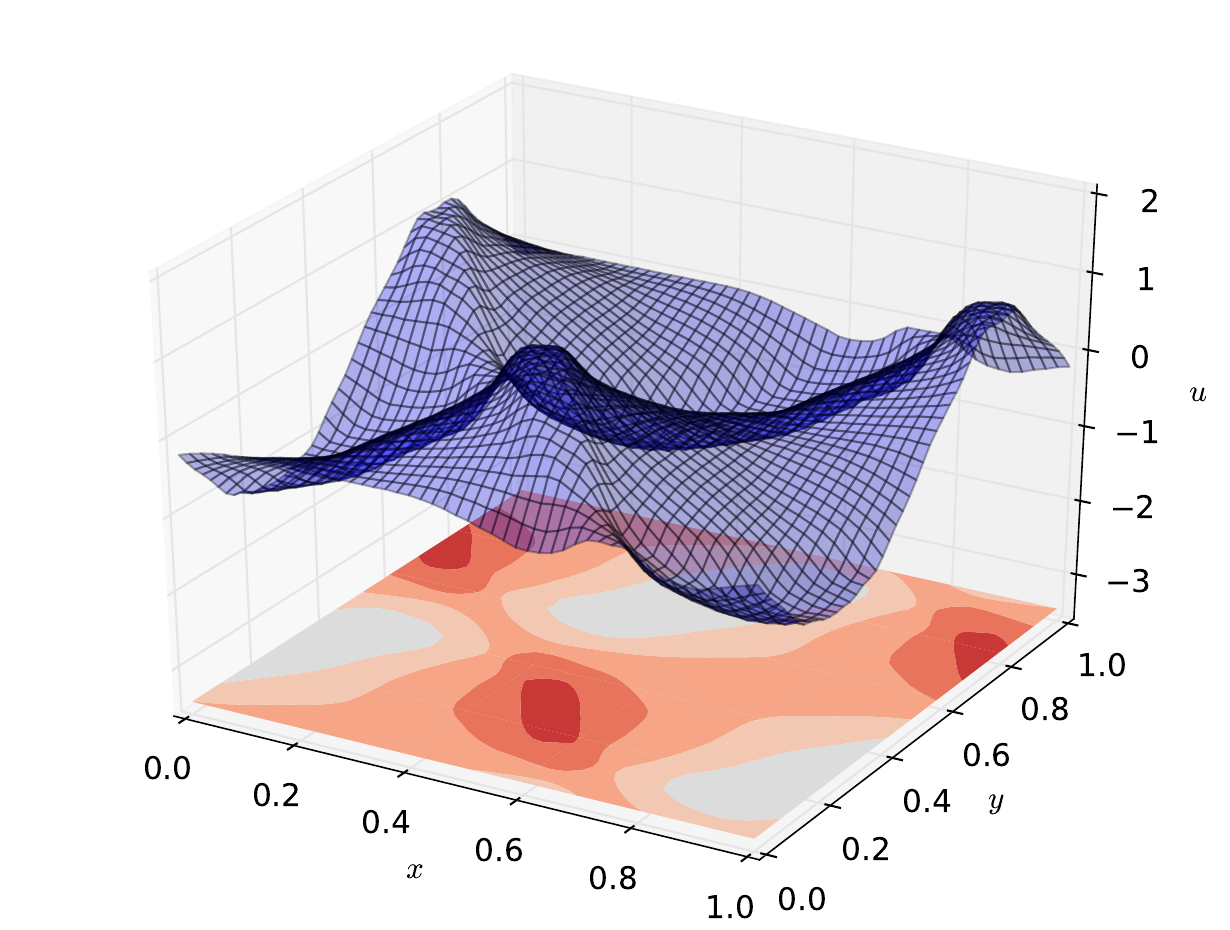}
    \caption{$t = 1$}
  \end{subfigure}
  \begin{subfigure}[b]{0.5\textwidth}
    \includegraphics[width=\textwidth]{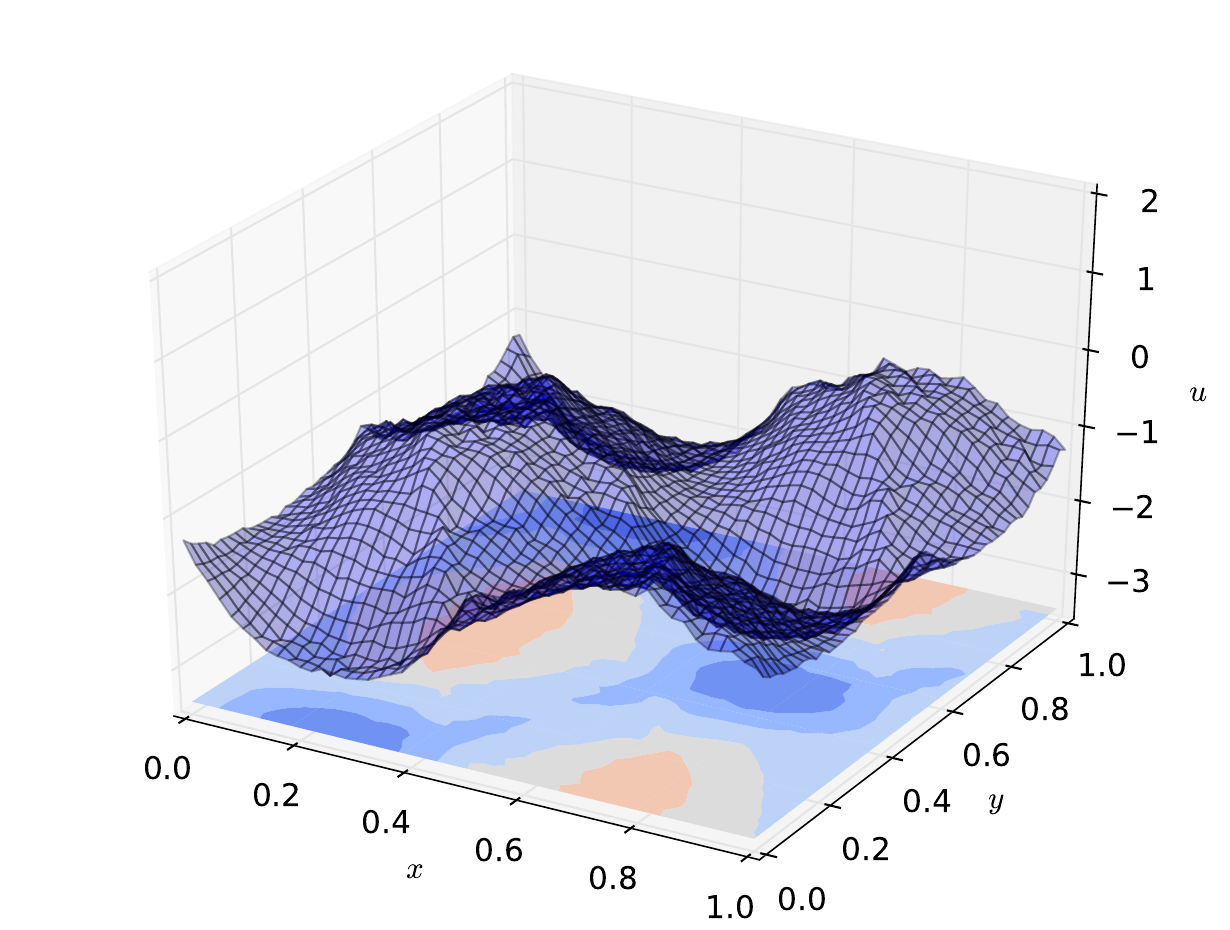}
    \caption{$t = 2$}
  \end{subfigure}
  \caption{Conservative DG scheme}
  \end{subfigure}
  \begin{subfigure}[b]{0.90\textwidth}
  \begin{subfigure}[b]{0.5\textwidth}
    \includegraphics[width=\textwidth]{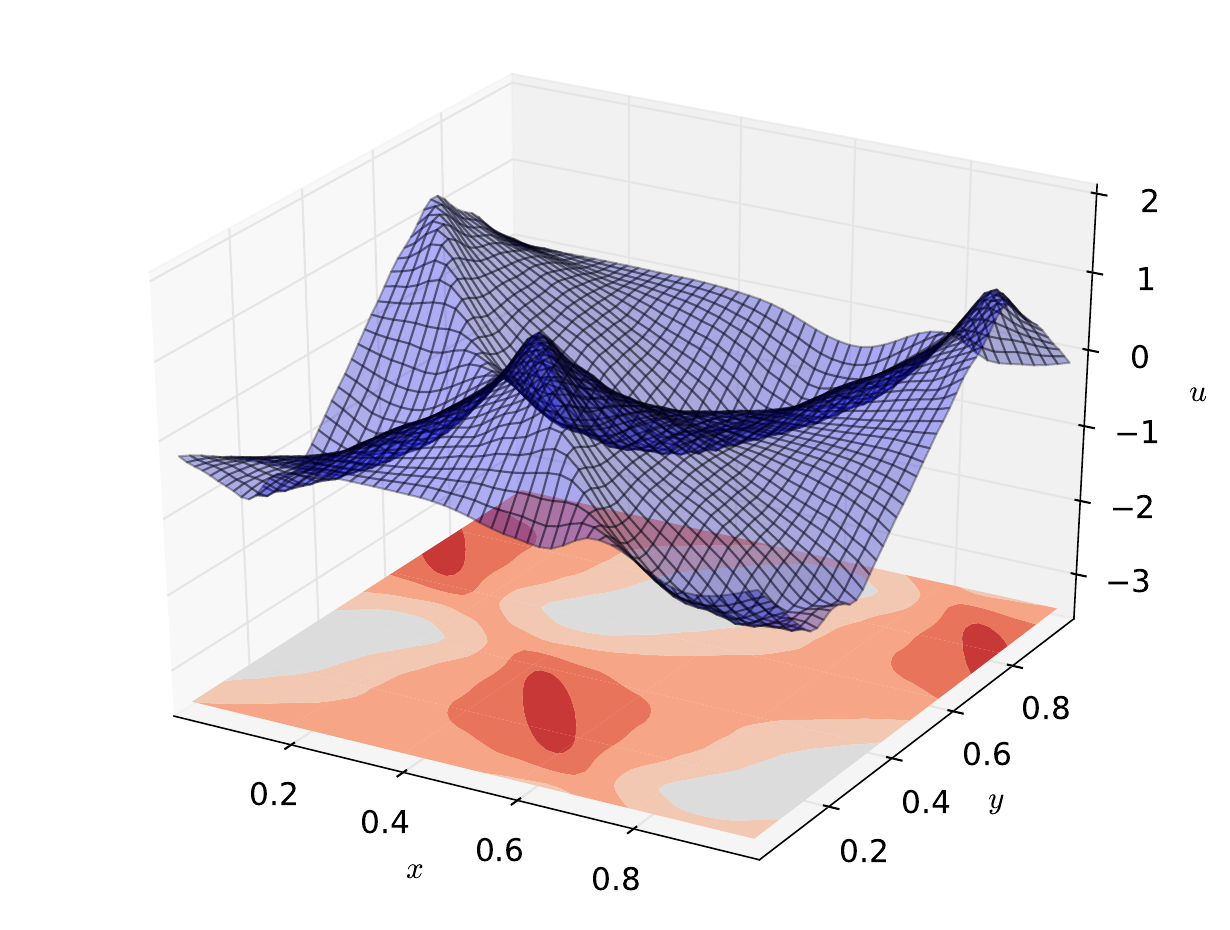}
    \caption{$t = 1$}
  \end{subfigure}
  \begin{subfigure}[b]{0.5\textwidth}
    \includegraphics[width=\textwidth]{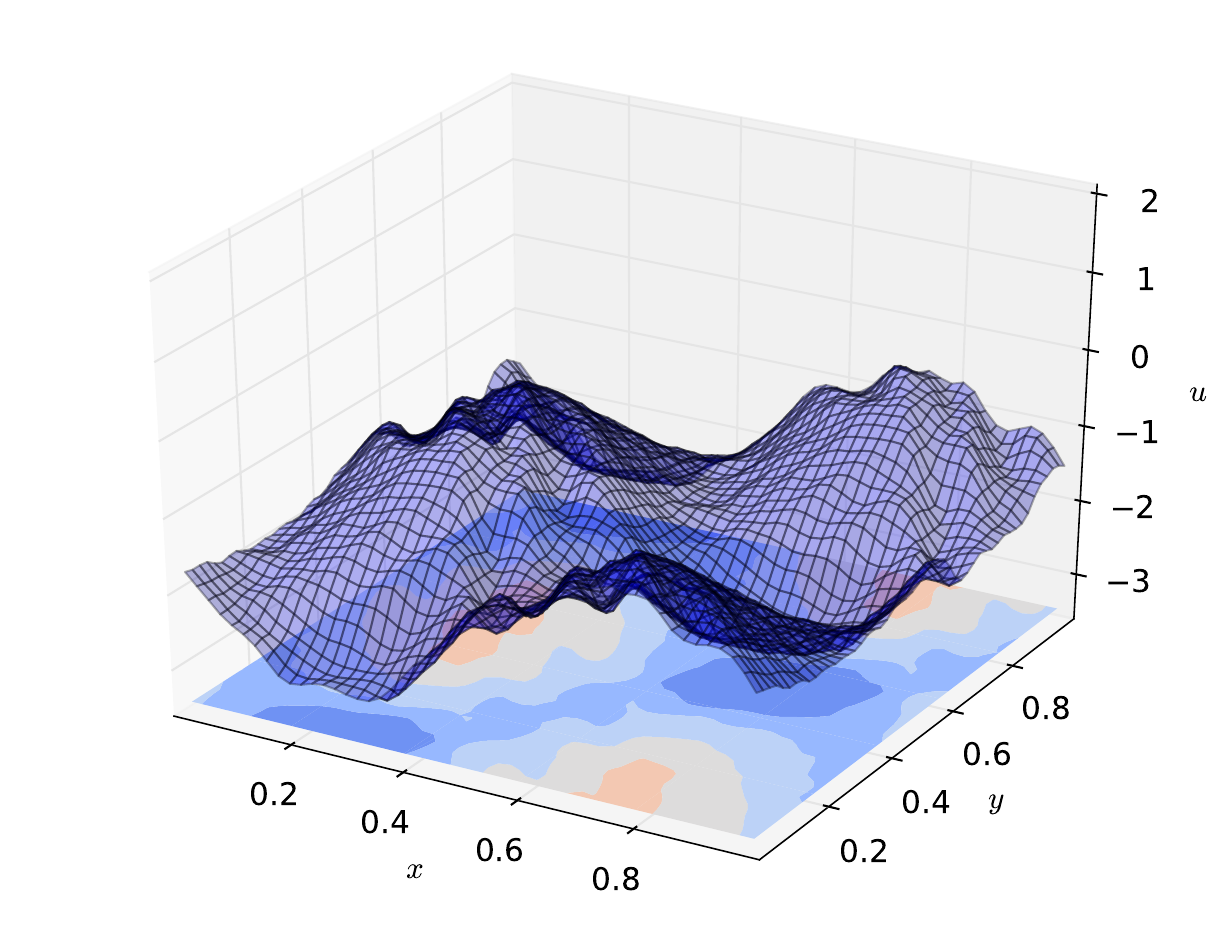}
    \caption{$t = 2$}
  \end{subfigure}
  \caption{Hamiltonian scheme}
  \end{subfigure}
  \caption{The numerical solution at left: $t=1$ and right: $t = 2$ 
    of the initial value problem
    \eqref{eq:main_2d} with initial data \eqref{eq:sineData} using the 
    conservative and dissipative piecewise quadratic 
    schemes ($s = 3$) with $N = 64$ cells. The bottom row shows the
    numerical solution using the Hamiltonian scheme.
    The physical parameters were $\alpha = 1.5$ and $\beta = 0.5$.}
  \label{fig:sw_case}
\end{figure}

Figure \ref{fig:sw_case}  shows the numerical solution using both conservative
and dissipative piecewise quadratic schemes with $N = 64$ cells. For
comparison, a numerical solution was also computed using the Hamiltonian
scheme derived in Section \ref{sec:ham_scheme}. The physical parameters were, as 
before, $\alpha = 1.5$ and $\beta = 0.5$. For $t > 0$ the non-isotropic
elasticity of the director field deteriorates the initial standing wave and
the pattern becomes more complicated. At $t = 2$ the solution given by the
dissipative DG scheme is visibly more regular that the solutions given by the 
conservative schemes (DG and Hamiltonian).

\section{Summary}

Using the Discontinuous Galerkin framework we have derived arbitrarily
high-order numerical schemes for the 2D variational wave equation describing
the director field in a type of nematic liquid crystals. By design, these
schemes either conserve or dissipate the total mechanical energy of the
system. The energy conserving scheme is based on a centralized numerical flux,
while the dissipative scheme employs a dissipative flux combined with a shock
capturing operator.

We have performed extensive numerical experiments both to verify the performance of
the schemes and to investigate the behavior of solutions to the variational
wave equation. In particular:
\begin{itemize}
  \item The schemes converge to a high order of accuracy for smooth solutions.
  \item The high-order schemes outperform low-order scheme in terms of 
    error per CPU time.
  \item The energy respecting properties (proven at the semi-discrete level)
    also hold on the fully discrete level when using a high-order numerical
    integration in time.
  \item Experiments show that the solution can loose regularity in finite time
    even for smooth initial data.
  \item After loss of regularity, results indicate that the conservative and dissipative 
    schemes converge to different solutions as the grid is refined. 
\end{itemize}

To the best of our knowledge, this is the first systematic numerical study of
the 2D generalization of the nonlinear variational wave equation
\eqref{eq:main}. Indeed, the results here indicate that the mathematical
treatment of \eqref{eq:main_2d} might be as delicate as in the 1D case.


\begin{thebibliography}{10}

\bibitem{Aursand2014Preprint} P.~Aursand and U.~Koley.
  \newblock Local discontinuous Galerkin schemes for a
      Nonlinear variational wave equation
        modeling liquid crystals,
  \newblock Preprint 2014

%\bibitem{ali2009} G.~Ali and J.K.~Hunter.
%   \newblock Orientational waves in a director field with rotational inertia,
%   \newblock {\em Amer. Inst. of Math. Sci.}, 2(1): 1--6 (2009).


%\bibitem{Badia1} S.~Badia, F.~Guill\'{e}n. Gonz\'alez and J.~V. Guti\'errez--Santacreu. 
%  \newblock An overview on Numerical Analyses of Nematic Liquid Crystal flows,
%   \newblock {\em Arch. Comput. Methods Eng.}, 18: 285--313 (2011).
%
%
%\bibitem{Badia2} S.~Badia, F.~Guill\'en. Gonz\'alez and J.~V. 
%   Guti\'errez--Santacreu.  
%  \newblock Finite element approximation of nematic liquid crystal 
%  flows using a 
%  saddle-point structure,
%  \newblock {\em J. Comput. Phys.}, 230(4): 1686--1706 (2011).


%\bibitem{Bartels} S.~Bartels and A.~Prohl.  \newblock
%  Constraint preserving implicit finite element discretization of harmonic map heat 
%  flow into spheres,
%  \newblock {\em Math. Comput.}, 76: 1847--1859 (2007).


\bibitem{barth} T. J. ~Barth.  \newblock
  Numerical methods for gas-dynamics systems on unstructured meshes. 
  In: \emph{An introduction to recent developments in theory and numerics of conservation laws}
   \newblock {\em Lecture notes in computational science and engineering}, vol(5), Springer, Berlin. Eds: D Kroner, M. Ohlberger, and C. Rohde, 1999.


%\bibitem{Becker} R.~Becker, X.~Feng and A.~Prohl.  \newblock
% Finite element approximations of the Ericksen--Leslie model for nematic liquid crystal 
% flow, 
% \newblock {\em SIAM J. Numer. Anal.}, 46(4): 1704--1731 (2008).


\bibitem{Beres} H.~Berestycki, J.~M. Coron and I.~Ekeland.  \newblock
  Variational Methods, Progress in nonlinear differential equations
  and their applications, \newblock {\em Vol 4}, Birkh\"{a}user,
  Boston, 1990.

\bibitem{bressan} A. ~Bressan and Y. ~Zheng.  \newblock Conservative
  solutions to a nonlinear variational wave equation, \newblock {\em
    Commun. Math. Phys.}, 266: 471--497 (2006).


\bibitem{chavent} G. ~Chavent and B. ~Cockburn.  \newblock
  The local projection $p^0 p^1$-discontinuous Galerkin finite element methods for 
  scalar conservation law,
  \newblock {\em Math. Model. Numer. Anal.}, 23: 565--592 (1989).



%\bibitem{chris} D. ~Christodoulou and A. ~Tahvildar-Zadeh.  \newblock
%  On the regularity of Spherically symmetric wave maps, 
%  \newblock {\em
%    Comm. Pure Appl. Math.}, 46: 1041--1091 (1993).


\bibitem{cockburnlinshu} S. Y. ~Cockburn, B. ~Lin and C. W. ~Shu.  \newblock
TVB Runge-Kutta local projection discontinuous Galerkin finite element methods for conservation laws III: one dimensional systems,
\newblock
  {\em J. Comput. Phys.}, 84: 90--113 (1989).


\bibitem{coron} J. ~Coron, J. ~Ghidaglia and F. ~H\'elein.  \newblock
  {\em Nematics}, Kluwer Academic Publishers, Dordrecht, 1991.

%\bibitem{diperna} R. ~J. Diperna and A. ~Majda.  \newblock
%  Oscillations and Concentrations in weak solutions of the
%  incompressible fluid equations, \newblock {\em Comm. Math. Phys},
%  108: 667--689 (1987).


\bibitem{ericksen} J. ~L. Ericksen and D. ~Kinderlehrer.  \newblock
  Theory and application of Liquid Crystals, \newblock {\em IMA
    Volumes in Mathematics and its Applications}, Vol 5, Springer
  Verlag, New York, 1987.

\bibitem{Gang1987} X. Gang, S. Chang-Qing, and L. Lei
  \newblock Perturbed solutions in nematic liquid crystals under 
  time-dependent shear.
  \newblock {\em Phys. Rew. A}, 36(1): 277--284 (1987).

\bibitem{glassey} R. ~T. Glassey.  \newblock Finite-time blow-up for
  solutions of nonlinear wave equations, 
  \newblock {\em Math. Z.}, 177: 1761--1794 (1981).

\bibitem{ghz1997} R. ~Glassey, J. ~Hunter, and Y. ~Zheng.  \newblock
  Singularities and Oscillations in a nonlinear variational wave
  equation. In: J.~Rauch and M.~Taylor, editors, \newblock {\em
    Singularities and Oscillations}, Volume 91 of the IMA volumes in
  Mathematics and its Applications, pages 37--60. Springer, New York,
  1997.

\bibitem{ghz1996} R. ~T. Glassey, J. ~K. Hunter and Yuxi. ~Zheng.
  \newblock Singularities of a variational wave equation, \newblock
  {\em J. Diff. Eq.}, 129: 49--78 (1996).

%\bibitem{gottliebetal} S.~Gottlieb, C.-W.~Shu and
%  E.~Tadmor. \newblock Strong stability preserving high-order time
%  discretization methods, \newblock {\em SIAM Review}, 43(1): 89-112
%  (2001).

\bibitem{hill} T. R. ~Hill and W. H. ~Reed.  \newblock Triangular mesh methods for neutron transport equation, \newblock
  {\em Tech. Rep. LA-UR-73-479.}, Los Alamos Scientific Laboratory, 1973.

\bibitem{hiltebrand} A. ~Hiltebrand and S. ~Mishra.  \newblock 
  Entropy stable shock capturing space--time discontinuous {G}alerkin schemes
  for systems of conservation laws, \newblock
  {\em Numer.~Math.} 126(1): 103--151 (2014).

\bibitem{holden} H. ~Holden and X. ~Raynaud.  \newblock Global
  semigroup for the nonlinear variational wave equation, \newblock
  {\em Arch. Rat. Mech. Anal.}, 201(3): 871--964 (2011).


\bibitem{hkr2009} H. ~Holden, K. ~H. Karlsen, and N. ~H. Risebro.
  \newblock A convergent finite-difference method for a nonlinear
  variational wave equation, \newblock {\em IMA. J. Numer. Anal.},
  29(3): 539--572 (2009).

\bibitem{hs1991} J. ~K. Hunter and R. ~A. Saxton.
  \newblock Dynamics of director fields, 
  \newblock {\em SIAM J. Appl. Math.},
  51: 1498--1521 (1991).

\bibitem{Johnson1990} C. Johnson, P. Hansbo and A. Szepessy,
  \newblock On the convergence of shock capturing streamline diffusion
  methods for hyperbolic conservation laws,
  \newblock {\em Math. Comput.},
  54(189): 107--129 (1990).

\bibitem{Kapustina2004} O. A. Kapustina.
  \newblock Liquid crystal acoustics: A modern view of the problem.
  \newblock {\em Crystallogr. Rep.} 49(4): 680--692 (2004)  

\bibitem{koley} U. ~Koley, S. ~Mishra, N. ~H. Risebro, and F. ~Weber.
  \newblock Robust finite-difference schemes for a nonlinear
  variational wave equation modeling liquid crystals, \newblock {\em Submitted}.

\bibitem{Leslie1979} F. ~M. Leslie.  \newblock 
Theory of flow phenomena in liquid crystals,
\newblock {\em Liquid Crystals}, 4, 1--81 (1979).

%  \bibitem{Lin} F. ~H. Lin and C.~Liu.  \newblock 
%  Non-parabolic dissipative systems modelling the flow of liquid crystals,
%  \newblock {\em Comm. Pure Appl. Math.}, 48, 501--537 (1995).

%\bibitem{Liu} F. ~H. Lin and C.~Liu. \newblock 
%Existence of solutions for the Ericksen-Leslie system,
% \newblock {\em 
% Arch. Ration. Mech. Anal.}, 154, 135--156 (2000).

\bibitem{Lut} H.~A. Luther and H.~P. Konen. \newblock 
 Some fifth-order classical Runge--Kutta formulas
 \newblock {\em SIAM Review}, 7(4): 551--558 (1965).

\bibitem{saxton} R. ~A. Saxton.  \newblock Dynamic instability of the
  liquid crystal director, \newblock {\em Contemporary Mathematics Vol
    100}, Current Progress in Hyperbolic Systems, pages 325--330,
  ed. W. B. Lindquist, AMS, Providence, 1989.

%\bibitem{shatah} J. ~Shatah.  \newblock Weak solutions and development
%  of singularities in the SU(2) $\sigma$-model, \newblock {\em
%    Comm. Pure Appl. Math.}, 41: 459--469 (1988).
%
%\bibitem{shatah1} J. ~Shatah and A. ~Tahvildar-Zadeh.  \newblock
%  Regularity of harmonic maps from Minkowski space into rotationally
%  symmetric manifolds, \newblock {\em Comm. Pure Appl. Math},
%  45: 947--971 (1992).

\bibitem{shu} C.-W. Shu. \newblock
    Different formulations of the discontinuous {G}alerkin method for 
    the viscous terms, \newblock In: 
    {\em Conference in Honor of Professor H.-C.
    Huang on the occasion of his retirement}, Science Press, 14--45, 2000.

\bibitem{Stewart2004} I. ~W. Stewart.  \newblock
  {\em The Static and Dynamic Continuum theory of liquid crystals: a mathematical
  introduction}, \newblock CRC Press, Boca Raton, 2004.

\bibitem{vanDoorn1975} C. Z. van Doorn.
  \newblock Dynamic behavior of twisted nematic liquidcrystal
  layers in switched fields.
  \newblock {\em J. Appl. Phys.}, 46: 3738--3745 (1975).

\bibitem{Vladimirov2007} V. A. Vladimirov and M. Y. Zhukov.
  \newblock Vibrational freedericksz transition in liquid crystals.
  \newblock {\em Phys. Rev. E}, 76:031706 (2007).

\bibitem{Yun1973} C. K. Yun.
  \newblock Inertial coefficient of liquid crystals: A proposal for its
  measurements.
  \newblock {\em Phys. Lett. A}, 45(2): 119--120 (1973).

\bibitem{zhang1} P. ~Zhang and Y. ~Zheng.  \newblock On oscillations
  of an asymptotic equation of a nonlinear variational wave equation,
  \newblock {\em Asymptot. Anal.}, 18(3): 307--327 (1998).

\bibitem{zhang2} P. ~Zhang and Y. ~Zheng.  \newblock Singular and
  rarefactive solutions to a nonlinear variational wave equation,
  \newblock {\em Chin. Ann. Math.}, 22: 159--170 (2001).

\bibitem{zhang3} P. ~Zhang and Y. ~Zheng.  \newblock Rarefactive
  solutions to a nonlinear variational wave equation of liquid
  crystals, \newblock {\em Commun. Partial Differ. Equ.}, 26: 381--419 (2001).

\bibitem{zhang4} P. ~Zhang and Y. ~Zheng.  \newblock Weak solutions to
  a nonlinear variational wave equation, \newblock {\em
    Arch. Rat. Mech. Anal.}, 166: 303--319 (2003).

\bibitem{zhang5} P. ~Zhang and Y. ~Zheng.  \newblock Weak solutions to
  a nonlinear variational wave equation with general data, \newblock
  {\em Ann. Inst. H. Poincar\'e Anal. Non Lin\'eaire}, 22: 207--226 
  (2005).

\bibitem{zhang6} P. ~Zhang and Y. ~Zheng.  \newblock On the global
  weak solutions to a nonlinear variational wave equation, \newblock
  {\em Handbook of Differential Equations. Evolutionary Equations},
  ed. C. M. Dafermos and E. Feireisl, vol. 2,
  pages 561--648, Elsevier, 2006.

\end{thebibliography}
\end{document}